\documentclass[11pt]{article}
\newcommand{\indic}[1]{\mathbf{1}_{\{#1\}}}

\usepackage{epsfig,latexsym,amssymb, subfigure}
\usepackage{amsmath} \usepackage{amsthm} \usepackage{amsfonts}

\newcommand{\eq}{\begin{equation}}
\newcommand{\en}{\end{equation}}
\newcommand{\re}[1]{\mbox{(\ref{#1})}}

\newtheorem{Theorem}{Theorem}
\newtheorem{theorem}[Theorem]{Theorem}
\newtheorem{lemma}[Theorem]{Lemma}

\newtheorem{corollary}[Theorem]{Corollary}
\newtheorem{construction}[Theorem]{Construction}
\newtheorem{proposition}[Theorem]{Proposition}
\newtheorem{example}{Example}
\newtheorem{exercise}[Theorem]{Exercise}
\newtheorem{defn}[Theorem]{Definition}

\newtheorem{question}[Theorem]{Question}
\newtheorem{conjecture}[Theorem]{Conjecture}
\newtheorem{condition}[Theorem]{Condition}
\newtheorem{remark}[Theorem]{Remark}
\newtheorem{problem}{Problem}

\def\endpf{\hfill $\Box$ \vskip .25in}

\newfont{\msbm}{msbm10 at 12pt}
\newfont{\eusb}{eusb10}
\newfont{\eusm}{eusm10}
\newfont{\eurb}{eurb10}
\newfont{\eurm}{eurm10}
\newfont{\eufb}{eufb10}
\newfont{\eufm}{eufm10}

\newcommand {\PR} {\mathbb{P}}

\newcommand{\te}{\rightarrow}
\newcommand{\ed}{\mbox{$ \ \stackrel{d}{=}$ }}

\newcommand{\comment}[1]{}
\newcommand{\noshowcomment}{\renewcommand{\comment}[1]{}}

\newenvironment{exm}[1]{\begin{example}\protect\label{#1}}{\end{example}}

\newenvironment{prp}[1]{\begin{proposition}\protect\label{#1}}{\end{proposition}}

\newenvironment{pbm}[1]{\begin{problem}\protect\label{#1}}{\end{problem}}

\newcommand{\Pnn}{{\cal P}_{[n]}}
\newcommand{\Pnk}{{\cal P}_{[n,k]}}
\newcommand{\Bnk}{{B_{n,k}}}
\newcommand{\Pn}{{\cal P}_{n}}
\newcommand{\TT}{\mathcal{T}}
\newcommand{\VV}{\mathbf{V}}
\newcommand{\TTn}{\mathcal{T}_n}
\newcommand{\Pisk}[1]{ { \Pi_{#1} ^ * } }

\def\mn{\medskip\noindent}

\def\beq{\begin{equation}}
\def\eeq{\end{equation}}
\def\beqa{\begin{eqnarray}}
\def\eeqa{\end{eqnarray}}
\def\beqax{\begin{eqnarray*}}
\def\eeqax{\end{eqnarray*}}
\def\th{^{\text{th}}}

\newcommand{\la}{\lambda}
\newcommand{\giv}{\left| \right.}
\noshowcomment
\begin{document}
\title{Gibbs distributions for random partitions generated by a
fragmentation process }
\author{Nathana\"el Berestycki and Jim Pitman\thanks{Research supported in part by N.S.F.
Grant DMS-0405779}
\\
\\
University of British Columbia, \\
and Department of Statistics, U.C. Berkeley }
\date{July 24, 2006}
\maketitle
\begin{abstract}
In this paper we study random partitions of $\{1,\ldots,n\}$ where
every cluster of size $j$ can be in any of $w_j$ possible internal
states. The Gibbs $(n,k,w)$ distribution is obtained by sampling
uniformly among such partitions with $k$ clusters. We provide
conditions on the weight sequence $w$ allowing construction of a
partition valued random process where at step $k$ the state has
the Gibbs $(n,k,w)$ distribution, so the partition is subject to
irreversible fragmentation as time evolves. For a particular
one-parameter family of weight sequences $w_j$, the time-reversed
process is the discrete Marcus-Lushnikov coalescent process with
affine collision rate $K_{i,j}=a+b(i+j)$ for some real numbers $a$
and $b$. Under further restrictions on $a$ and $b$, the
fragmentation process can be realized by conditioning a
Galton-Watson tree with suitable offspring distribution to have
$n$ nodes, and cutting the edges of this tree by random sampling
of edges without replacement, to partition the tree into a
collection of subtrees. Suitable offspring distributions include
the binomial, negative binomial and Poisson distributions.
\end{abstract}


\vfill \noindent \textbf{Keywords} Fragmentation processes, Gibbs
distributions, Marcus-Lushnikov processes, Gould convolution
identities.

\newpage

\section{Introduction}

Gibbs models for random partitions generated by random processes
of coagulation and fragmentation have been widely studied
(\cite{whittle1}, \cite{whittle2}, \cite{whittle3}, \cite{kelly}).
They typically arise as equilibrium distributions of
time-reversible processes of coagulation and fragmentation (see
for instance \cite{durrett-granovsky-gueron}, and \cite{jb} for
general results about exchangeable fragmentation-coalescence
processes in equilibrium). There is a much smaller literature in
which Gibbs models are derived from an irreversible Markovian
coagulation process \cite{hses85}. This paper presents a Gibbs
model for an irreversible Markovian fragmentation process. While
Gibbs models for physical processes of fragmentation have been
treated before, such models typically allow the possibility of
both fragmentation and coagulation at the microscopic level,
resulting in a Gibbs equilibrium at any given time. However as
time evolves, this equilibrium is moving towards a more fragmented
state, so using the language of thermodynamics this equilibrium
should be seen as a quasistatic equilibrium of an adiabatic
process. The point here is to provide a rigorous Markovian model
of irreversible fragmentation at the microscopic level with no
possibility of coagulation allowed.

\medskip The simplest way to describe the process treated here is to
specify its time reversal. This is the Marcus-Lushnikov coalescent
process with collision rate kernel $K_{i,j} = a + b (i + j)$ for
some constants $a$ and $b$, where $K_{i,j}$ represents rate of
collisions between clusters of $i$ particles and clusters of $j$
particles. This model was solved by Hendriks et al. \cite{hses85},
who showed that the distribution at time $t$ in such a coalescent
process started from a monodisperse initial condition is a mixture
of microcanonical Gibbs distributions with mixing coefficients
depending on $t$. Here, we derive what turns out to be essentially
the same model, modulo time reversal and a formulation in discrete
rather than continuous time, but from a different set of
assumptions describing the evolution of the process with time
running in the direction of fragmentation. The probabilistic link
between the two sets of assumptions is a time reversal calculation
using Bayes rule. The most interesting feature of this calculation
is that starting from a natural recursive assumption for the
fragmentation process in terms of Gibbs distributions, there is
only one-parameter family of possible solutions to the problem,
with the parameter corresponding to the ratio of the two
parameters $a$ and $b$ in the collision rate kernel of the
reversed time process. Another interesting feature is that for
some but not all $a$ and $b$, the fragmentation process for
partitions of a set of size $n$ can be realized by conditioning a
Galton-Watson process with suitable offspring distribution to have
a family tree of size $n$, then cutting the edges of this tree by
a process of random sampling without replacement. The suitable
offspring distributions include the binomial, negative binomial
and Poisson distributions.

\subsection{Canonical and microcanonical Gibbs distribution}

\medskip Typically, the state of a coagulation/fragmentation process
is represented by a random partition of $n$, that is a random
variable with values in the set $\Pn$ of all partitions of $n$. In
later sections of this paper the state of the process will be
represented rather as a random partition of the set $\{1,2,
\ldots, n \}$, as this device simplifies a number of calculations.
But the rest of this introduction follows the more common
convention of working with the set $\Pn$ of partitions of the
integer $n$. Let
\begin{equation}
\label{lam} \la = 1^{c_1} 2^{c_2} \cdots n^{c_n}
\end{equation}
denote a typical partition of $n$. Regarding the state of the
system as a partition of $n$ particles into clusters of various
sizes, the state $\la$ in (\ref{lam}) indicates that there are
$c_j$ clusters of size $j$ for each $1 \le j \le n$. Note that
$\sum_j j c_j = n $, the total number of particles. The total
number of clusters is $k:= \sum_j c_j$. The numbers $c_1, c_2,
\ldots, c_j \ldots$ may be called numbers of {\em monomers,
dimers, } $\ldots$ {\em j-mers}, or numbers of {\em singletons,
doubletons, } $\ldots$ {\em j-tons}. The Gibbs model most commonly
derived from equilibrium considerations is the {\em canonical
Gibbs distribution on partitions of $n$ with weight sequence
$(w_j)$} defined by
\begin{equation} \label{cangibbs} P(\la \giv
n ; w_1, w_2, \ldots, w_n) = { n! \over Y_n } \, \prod_{ i = 1
}^n\, {1 \over c_i! } \, \left( {w_i \over i !} \right) ^{c_i}
\end{equation}
where
\begin{equation} Y_n = Y_n(w_1, w_2, \ldots, w_n) \end{equation}
is a normalization constant. This polynomial in the first $n$
weights $w_1, w_2, \ldots, w_n$ is known in the combinatorics
literature as the {\em complete Bell } (or {\em exponential}) {\em
polynomial} \cite{comtet}. In the physics literature the Gibbs
formula (\ref{cangibbs}) is commonly written in terms of $x_i =
w_i/i!$ instead of $w_i$, and the polynomial
\begin{equation} Z_n(x_1,x_2, \ldots, w_n):= n! Y_n (1!  x_1, 2! x_2 ,
\ldots n! x_n)
\end{equation} is called the {\em canonical partition
function}. For textbook treatments of such models, and references
to earlier work see \cite{saintflour}. Typically, the canonical
Gibbs distribution (\ref{cangibbs}) is derived either from
thermodynamic considerations, or from a set of detailed balance
equations corresponding to a reversible equilibrium between
processes of fragmentation and coagulation. In the latter case the
canonical Gibbs distribution is represented as the equilibrium
distribution of a time-reversible Markov chain with state space
$\Pn$. For related models, see \cite{abt}, and Vershik
\cite{vershik}, who also considers a variation of Gibbs
distributions in the context of quantum statistical physics, where
he derives asymptotics for the limiting shape of a Gibbs partition
associated with the Bose-Einstein statistics by considering a
certain variational problem. See also \cite{gnedin-pitman} where
partitions are subject to a natural additional constraint of
consistency corresponding to infinite exchangeability.

\medskip Conditioning a canonical Gibbs distribution on the number
of clusters $k$ yields a corresponding {\em microcanonical Gibbs
distribution} for each $1\le k \le n$. This distribution assigns
to the partition $\la$ displayed in (\ref{lam}) the probability
\begin{equation} \label{micgibs} P(\la \giv n, k; w_1,
w_2, \ldots,w_n) = {n! \over \Bnk } \prod_{ i = 1 } ^ n \, {1
\over c_i! } \, \left( {w_i \over i !} \right) ^{c_i}
\end{equation}
where $\Bnk = \Bnk (w_1, w_2, \ldots )$ is a {\em partial Bell (or
exponential) polynomial}, and
\begin{equation} Z_{n,k} (x_1,x_2,
\ldots,x_n):= n! B_{n,k} (1! x_1, 2! x_2 , \ldots, n!x_n )
\end{equation}
is known as a {\em microcanonical partition function}. A great
many expressions, representations and recursions for these
polynomials $B_{n,k}$ and $Z_{n,k}$ are known \cite{comtet,
saintflour}. These formulae are useful whenever the weight
sequence $(w_j)$ is such that the associated polynomials admit an
explicit formula as functions of $n$ and $k$, or can be suitably
approximated (see e.g. \cite{holst}). Some of these results are
reviewed in \cite{saintflour}. The class of weight sequences
$(w_j)$ for which the microcanonical Gibbs model is ``solvable'',
meaning there is an explicit formula for the $\Bnk$, is quite
large.

\medskip
In the study of irreversible partition-valued processes, it has
been found in several cases (discussed below) that the
distribution of the process at time $t$ is a probabilistic mixture
over $k$ of microcanonical Gibbs distributions, that is to say a
probability distribution of the form
\begin{equation}
\label{gibbsform} P(\la) = \sum_{k=1}^n q_{n,k} P(\la \giv n, k;
w_1, w_2, \ldots, w_n)
\end{equation} where $q_{n,k}$ represents the probability that
$\la$ has $k$ components, so $q_{n,k} \ge 0 $ and $\sum_{k=1}^n
q_{n,k} = 1$, and both $q_{n,k}$ and the weight sequence $w_i$ may
be functions of $t$. We call any distribution of the form
\re{gibbsform} a {\em Gibbs distribution with weights $(w_j)$},
thereby including both canonical and microcanonical Gibbs
distributions. For example, Lushnikov \cite{lushnikov} showed that
the coalescent model with monodisperse initial condition and
collision rates $K_{x,y}  = x f(y ) + y f(x)$ leads to such Gibbs
distributions. Hendriks et. al \cite{hses85} showed that this is
also the case for $K_{x,y}  =  a + b ( x + y )$ for constants $a$
and $b$.  See also \cite[Section 1.5]{saintflour} for an
interpretation of Gibbs distributions with integer weights in
terms of composite combinatorial structures.

\subsection{Organization and summary of the paper.}

The rest of this paper is organized as follows. Section
\ref{prelim} presents some background material, and introduces the
formalism of Gibbs distributions over partitions of the set $\{1,
2, \ldots, n \}$. Results for irreversible fragmentation processes
are presented in Section \ref{fragment}. Section \ref{main result}
presents the main result of the paper. This result states that,
under an additional set of assumptions (most notably, the
\emph{linear selection rule}), it is possible to construct a Gibbs
fragmentation process with weight sequence $(w_j)$ if and only if
$w_j=\prod_{m=2}^{j}(mc +jb)$ for some constants $b$ and $c$: in
this case it is shown that the time-reversal of the fragmentation
process is the discrete Marcus-Lushnikov coalesent with affine
coalescent rate: $K_{i,j}=a+b(i+j)$. Section 5 provides some
background material on generating functions and branching
processes which is needed for the evaluation of a particular Bell
polynomial. Section 6 then presents the proofs of the results of
Section 4. These results leave out an important case, which is
that of the sequence $w_j=(j-1)!$. In section \ref{continuous}, we
approach this problem from the angle of Kingman's coalescent
process and the Ewens sampling formula. In particular we construct
a continuous analogue of the desired process. However, we show
that the existence of this process with discrete time cannot be
obtained by taking the discrete-time chain embedded in the
continuous process, so that its existence remains an open
question.

\section{Preliminaries}
\label{prelim} Let $[n]$ denote the set $\{1, \ldots, n \}$. A
{\em partition of $[n]$} is an unordered collection of non-empty
disjoint subsets of $[n]$ whose union is $[n]$. A generic
partition of $[n]$ into $k$ sets (sometimes also referred to as
blocks or components of the partition) will be denoted
\begin{equation}
\label{pik} \pi_k = \{A_1, \ldots, A_k \}    \mbox{, where } 1 \le
k \le n
\end{equation}
and where blocks are numbered according e.g. to their least
element. Let $\Pnn$ denote the set of all partitions of $[n]$, and
let $\Pnk$ be the subset of $\Pnn$ comprising all partitions of
$[n]$ into $k$ components. Given a sequence of weights $(w_j, j =
1, 2, \cdots )$ define the {\em microcanonical Gibbs distribution
on $\Pnk$ with weights $(w_1, w_2, \ldots)$} to be the probability
distribution on $\Pnk$ which assigns to each partition $\pi_k$ as
in (\ref{pik}) the probability
\begin{equation}
\label{gibbsnk} p_{n,k}( \pi_k ; w_1, w_2, \ldots w_n) = {1 \over
\Bnk} \, \prod_{i=1}^k w_{\# A_i},
\end{equation}
where $\#A_i$ denotes the number of elements of $A_i$ and
\begin{equation} \label{belldef}
\Bnk := \Bnk (w_1, \ldots, w_{n} )  := \sum_{\pi_k \in
\mathcal{P}_{[n,k]}} \, \prod_{i=1}^k w_{n_i(\pi_k)},
\end{equation}
where for $\pi_k \in \mathcal{P}_{[n,k]}$, the $n_i (\pi_k)$ for
$1 \le i \le k$ are the sizes of the components of $\pi_k$ in some
arbitrary order. Throughout this paper we will use the notation
$p_{n,k}( \, \cdot \, ; w_1,\ldots,w_{n})$, or simply $p_{n,k}$ if
no confusion is possible, for the microcanonical Gibbs
distribution (\ref{gibbsnk}) determined by the sequence
$(w_1,\ldots, w_{n})$. Remark that as soon as $k\ge 2$, the
microcanonical Gibbs distribution $p_{n,k}$ actually only depends
on $(w_1,\ldots,w_{n-1})$. Given a partition $\pi$ of $[n]$, the
{\em corresponding partition $\la$ of $n$} is $\la := 1^{c_1}
2^{c_2} \cdots n^{c_n}$ as in (\ref{lam}) where $c_i$ is the
number of components of $\pi$ of size $i$. For each vector of
non-negative integer counts $(c_1, \ldots, c_n)$ with $\sum_i i
c_i = n$ the number of partitions $\pi$ of $[n]$ corresponding to
the partition $1^{c_1} 2^{c_2} \cdots n^{c_n}$ of $n$ is well
known to be
\begin{equation} \label{nfactor} { n ! \over \displaystyle \prod_{i=1}^n
c_i! (i! ) ^{c_i } } .
\end{equation}
The probability distribution on partitions of $n$ induced by the
microcanonical Gibbs distribution on $\Pnk$ with weights $(w_1,
w_2, \ldots)$ therefore identical to the microcanonical Gibbs
distribution on $\Pn$ with the same weights $(w_1, w_2, \ldots)$
as defined in (\ref{micgibs}), and there is the following standard
expression for $B_{n,k}$ \cite{comtet}:
\begin{equation}
\label{bnkdef} \Bnk = {n!} \sum_{ \la_k} \prod_{ i = 1 } ^ n \, {1
\over c_i! } \, \left( {w_i \over i !} \right) ^{c_i}
\end{equation}
where the sum is over all partitions $\la_k$ of $n$ into $k$
components, and $c_i = c_i(\la_k)$ is the number of components of
$\la_k$ of size $i$. Thus transferring from Gibbs distributions on
partitions of $n$ into $k$ components to Gibbs distributions on
partitions of the set $[n]$ into $k$ components is just a matter
of keeping track of the universal combinatorial factor
(\ref{nfactor}).

\subsection{Combinatorial Interpretation}

The following well-known interpretations provide both motivation
and intuition for the study of Gibbs distributions and Bell
polynomials. Suppose that $n$ particles labelled by elements of
the set $[n]$ are partitioned into {\em clusters} in such a way
that each particle belongs to a unique cluster. Formally, the
collection of clusters is represented by a partition of $[n]$.
Suppose further that each cluster of size $j$ can be in any one of
$w_j$ different {\em internal states} for some sequence of
non-negative integers $(w_j)$. Let the {\em configuration} of the
system of $n$ particles be the partition of the set of $n$
particles into clusters, together with the assignment of an
internal state to each cluster. For each partition $\pi$ of $[n]$
with $k$ components of sizes $n_1, \ldots, n_k$, there are
$\prod_{i=1}^k w_{n_i}$ different configurations with that
partition $\pi$. So $\Bnk (w_1, w_2 , \ldots)$ defined by
(\ref{belldef}) gives {\em the number of configurations with $k$
clusters}; the Gibbs distribution (\ref{gibbsnk}) with weight
sequence $(w_j)$ is {\em the distribution of the random partition
of $[n]$ if all configurations with $k$ clusters are equally
likely}, and formula (\ref{micgibs}) describes the corresponding
Gibbs distribution on partitions of $n$ induced by the same
hypothesis.

Many particular choices of $(w_j)$ have natural interpretations,
both combinatorial and physical. In particular, the following four
examples have been extensively studied. Many more combinatorial
examples are known where Gibbs distributions arise naturally from
an assumption of equally likely outcomes on a suitable
configuration space. Related problems of enumeration and
asymptotic distributions have been extensively studied
\cite{holst, kolchin, friedst, goh, renyi}.

\subsection{Some important examples}

We recall here some natural examples of Gibbs distributions for
particular sequences of weights $(w_j)$, and their combinatorial
interpretations, which motivate our work in following sections.

\begin{exm}{exm1}
Uniform distribution on partitions of $[n]$. {\em Take $w_j = 1$
for all $j$. Then a configuration is just a partition of $[n]$, so
that $B_{n,k} (1,1, \ldots )$ is the number of partitions of $[n]$
into $k$ components, known as a {\em Stirling number of the second
kind}. The microcanonical Gibbs model
$p_{n,k}$ corresponds to assuming that all partitions of $[n]$
into $k$ components are equally likely. }
\end{exm}

\begin{exm}{exm2}
Uniform distribution on permutations. {\em Suppose that the
internal state of a cluster $C$ of size $j$ is one of the $(j-1)!$
cyclic permutations of $C$. Then $w_j=(j-1)!$, and each
configuration corresponds to a permutation of $[n]$. Therefore
$B_{n,k} (0!,1!,2! \ldots )$ is the number of permutations of
$[n]$ with $k$ cycles, known as an {\em unsigned Stirling number
of the first kind}. The microcanonical Gibbs distribution
$p_{n,k}$
is the distribution on $\Pn$ induced by a permutation uniformly
chosen among all permutations with $k$ cycles.}
\end{exm}

\begin{exm}{exm3} Cutting a rooted random segment {\em
\cite{saintflour}. Suppose that the internal state of a cluster
$C$ of size $j$ is one of $j!$ linear orderings of the set $C$.
Identify each cluster as a directed graph in which there is a
directed edge from $a$ to $b$ if and only if $a$ is the immediate
predecessor of $b$ in the linear ordering. Call such a graph a
{\em rooted segment}. \begin{figure}[ht] \begin{center}
\epsfig{file=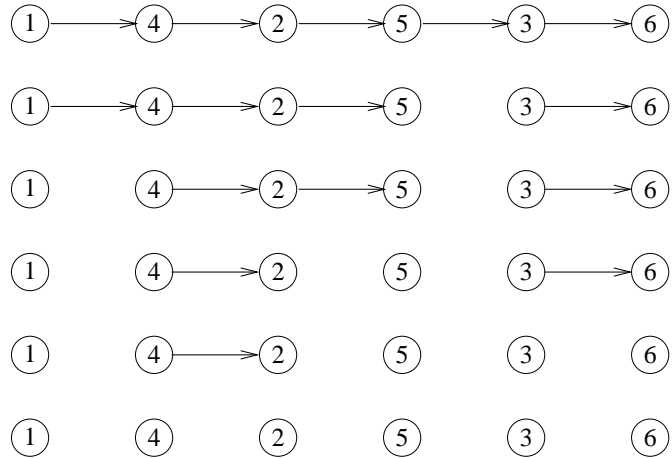, width=.6\textwidth,
height=.413\textwidth} \end{center} \caption{Cutting a rooted
random segment. } \label{figcutseg} \end{figure} Then $B_{n,k}
(1!,2!,3! \ldots )$ is the number of directed graphs whose
vertices are labelled by $[n]$ with $k$ such rooted segments as
its components. In the previous two examples, explicit formulae
for the $\Bnk$ are fairly complicated. But this time there is a
simple formula: \begin{equation} \label{lah} B_{n,k} (1!,2!,3!
\ldots ) = { n - 1 \choose k - 1 } \, {n! \over k ! }
\end{equation} is known as a {\em Lah number} \cite[p.
135]{comtet}. The Gibbs model in this instance is a variation of
Flory's model for a linear polymerization process \cite{flory}.
Another interpretation is provided by Kingman's coalescent
\cite{aldous, ki82co}. It is easily shown in this case that a
sequence of random partitions $(\Pi_k, 1 \le k \le n)$ such that
$\Pi_k$ has the microcanonical Gibbs distribution with $k$ blocks,
can be obtained as follows. Let $G_1$ be a uniformly distributed
random rooted segment labelled by $[n]$, and let $G_k$ be derived
from $G_1$ by deletion of a set of $k-1$ edges picked uniformly at
random from the set of $n-1$ edges of $G_1$, and let $\Pi_k$ be
the partition induced by the components of $G_k$. If the $n-1$
edges of $G_1$ are deleted sequentially, one by one, the random
sequence $(\Pi_1, \Pi_{2}, \ldots, \Pi_n )$ is a refining sequence
of random partitions such that $\Pi_k$ has the Gibbs
microcanonical distribution (\ref{gibbsnk}). This is illustrated
in figure \ref{figcutseg}. The time-reversed sequence $(\Pi_n,
\Pi_{n-1}, \ldots, \Pi_1 )$ is then governed by the rules of {\em
Kingman's coalescent}: conditionally given $\Pi_{k}$ with $k$
components, $\Pi_{k-1} $ is equally likely to be any one of the
${k \choose 2}$ different partitions of $[n]$ obtained by merging
two of the components of $\Pi_k$. Equivalently, the sequence
$(\Pi_1, \Pi_{2}, \ldots, \Pi_n )$ has uniform distribution over
the set ${\cal R}_n$ of all refining sequences of partitions of
$[n]$ such that the $k$th term of the sequence has $k$ components.
The consequent enumeration $\# {\cal R}_n = n! (n-1)!/ 2^{n-1}$
was obtained by Erd\"os et al \cite{erdgm75}. The fact that
$\Pi_k$ determined by this model has the microcanonical Gibbs
distribution with $k$ blocks and weight sequence $w_j = j!$ was
obtained by Bayewitz et. al. \cite{BYKS74} and Kingman
\cite{ki82co}.
}
\end{exm}

\begin{exm}{exm4}
Cutting a rooted random tree {\em \cite{saintflour}. Suppose the
internal state of a cluster $C$ of size $j$ is one of the
$j^{j-1}$ rooted trees labelled by $C$.
\begin{figure}[ht]
\begin{center}
\mbox{\subfigure{\epsfig{file=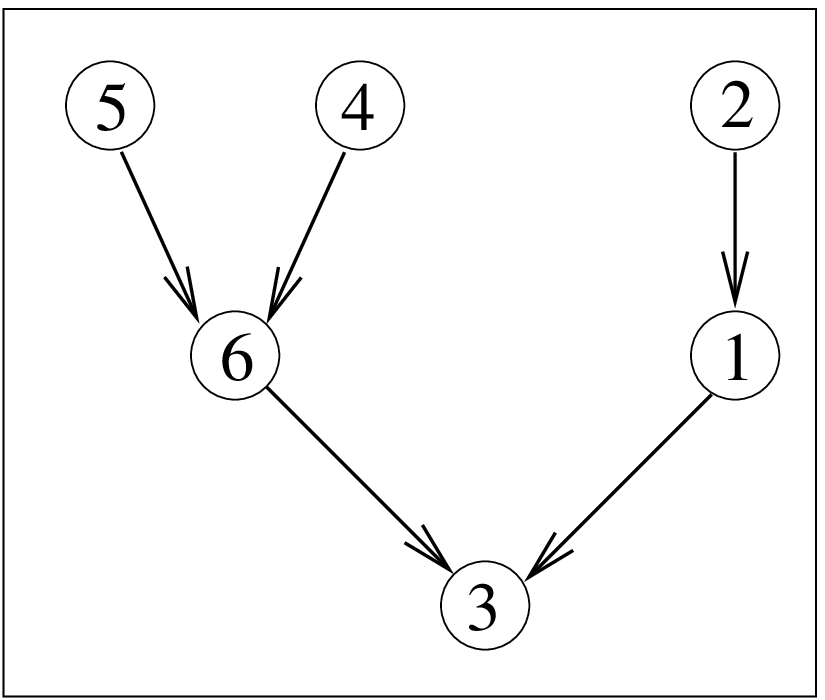, width=0.152\textwidth,
     height=0.128\textwidth}}\,
     \subfigure{\epsfig{file=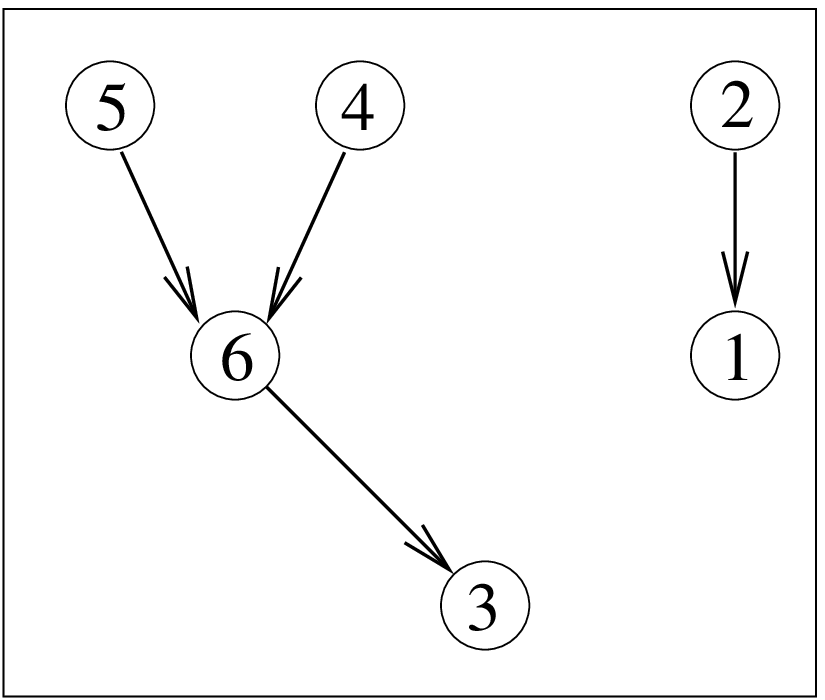, width=0.152\textwidth,
     height=0.128\textwidth}}\,
     \subfigure{\epsfig{file=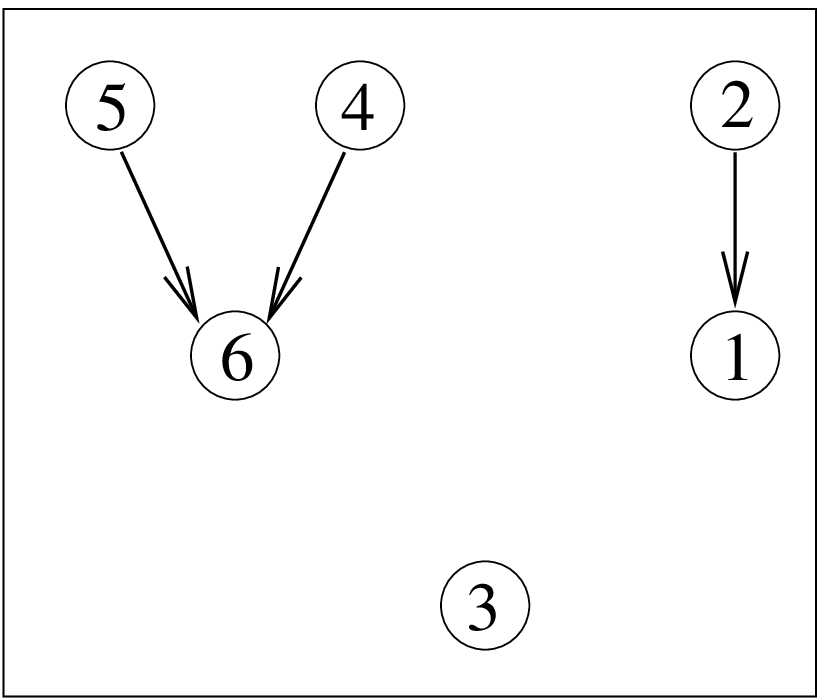, width=0.152\textwidth,
     height=0.128\textwidth}}\,
     \subfigure{\epsfig{file=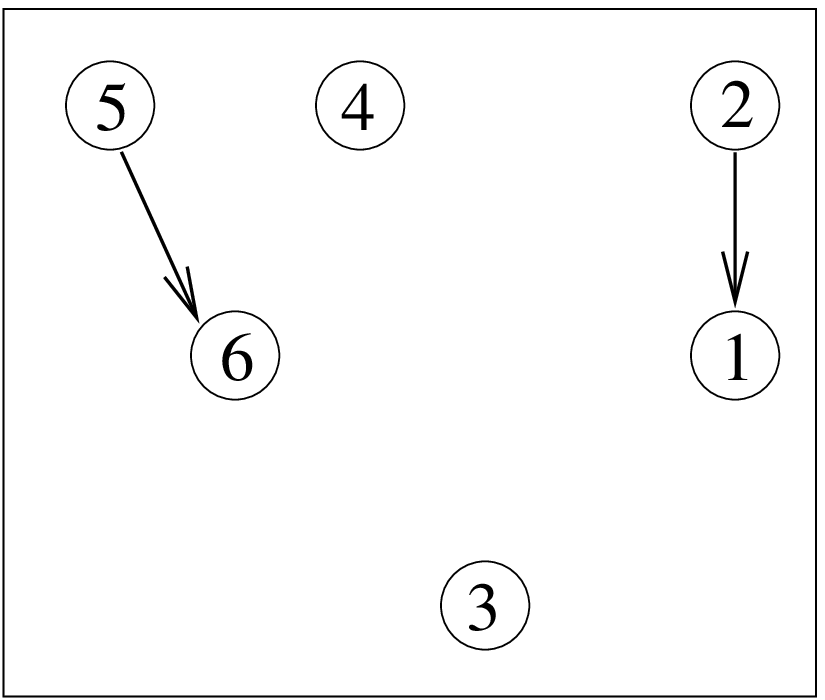, width=0.152\textwidth,
     height=0.128\textwidth}}\,
     \subfigure{\epsfig{file=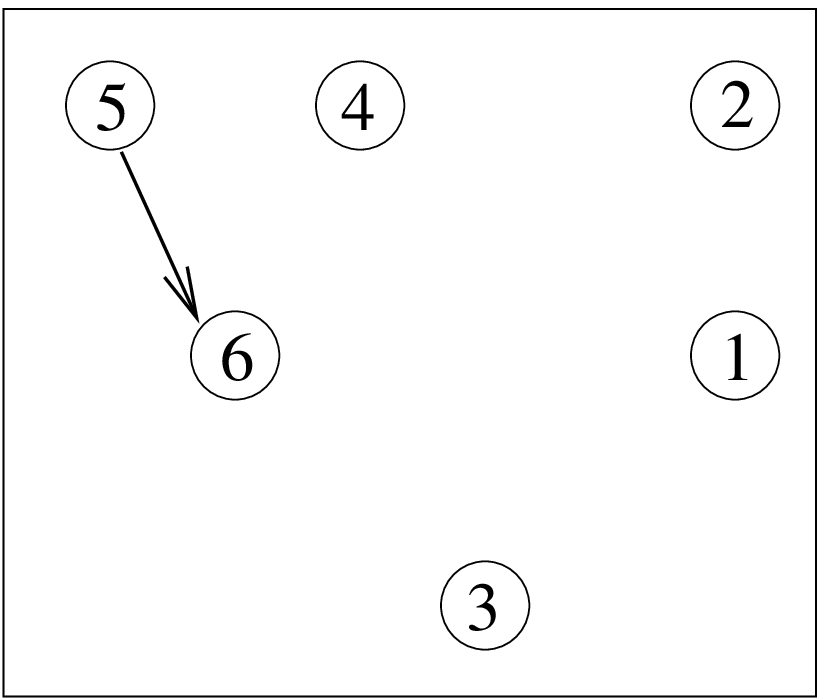, width=0.152\textwidth,
     height=0.128\textwidth}}\,
     \subfigure{\epsfig{file=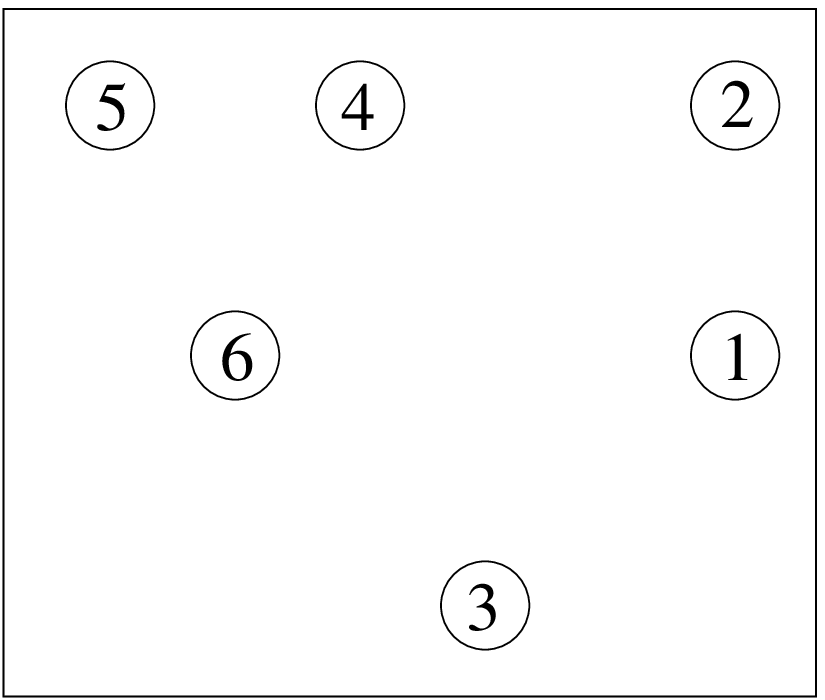, width=0.152\textwidth,
     height=0.128\textwidth}}}
\end{center}
\caption{Cutting a rooted random tree with 5 edges}
\label{figcuttree}
\end{figure}
Then $B_{n,k} ( 1^{1-1}, 2^{2-1}, 3^{3-1}, \ldots)$ is the number
of forests of $k$ rooted trees labelled $[n]$. This time again
there is a simple formula for $B_{n,k}$. As a consequence of
Cayley's enumeration of forests \cite{enumeration of
trees,randomforests}
\begin{equation} \label{cayley} B_{n,k} ( 1^{1-1}, 2^{2-1}, 3^{3-1},
\ldots) = { n - 1 \choose k - 1 } \, n^{n-k}
\end{equation}
The Gibbs model in this instance corresponds to assuming that all
forests of $k$ rooted trees are equally likely. This model turns
up naturally in the theory of random graphs and has been studied
and applied in several other contexts. The coalescent obtained by
reversing the process of deleting the edges at random is the
\emph{additive coalescent} as discussed in \cite{randomforests}.
}
\end{exm}

%

\section{Fragmentation Processes}
\label{fragment}

Recall that $\Pnn$ is the set of partitions of $[n]:= \{1, \ldots,
n\}$. Call a $\Pnn$-valued random process $(\Pi_t, t \in I)$, with
index set $I$ a subset of real numbers, a {\em fragmentation
process} if with probability one both
\begin{enumerate}
\item[(i)] for every pair of times $s$ and $t$ in $I$ with $s < t$
the partition $\Pi_t$ is a refinement of $\Pi_s$, and
 \item[(ii)]
for each $1 \le k \le n$ there is some $t \in I$ such that $\Pi_t$
has $k$ components.
\end{enumerate}

When a confusion is possible we will use the notation $\Pi(t)$ if
$I$ is a continuous interval and $\Pi_t$ in the case where $I$ is
a discrete subset of the real numbers. We emphasize that
throughout the paper, due to condition (ii), the fragmentation
processes we consider are binary. In other words, whenever a split
occurs, the split is a binary split in which one and only one
block of the partition splits in two, thereby incrementing the
number of components by 1. This condition also forces $\Pi(t)$ to
be the partition of $[n]$ with one component of size $n$ for all
sufficiently small $t\in I$, and to be the partition of $[n]$ into
$n$ singletons for all sufficiently large $t \in I$.

Given a sequence of numbers $(w_1,w_2,\ldots,w_{n-1})$, call
$(\Pi_t , t \in I)$ a {\em Gibbs fragmentation process with
weights} $(w_1, \ldots, w_{n-1})$ if for every $t \in I$ and $1
\le k \le n$, the conditional distribution of $\Pi_t$ given that
$\Pi_t$ has $k$ components is the microcanonical Gibbs
distribution $p_{n,k}$ on $\Pnk$ as defined by (\ref{gibbsnk}).
Note that if $(\Pi_{k}, k \in [n])$ is a Gibbs fragmentation
process, then the unconditional distribution of $\Pi_k$ is also
$p_{n,k}$, because condition (ii) implies that $\Pi_{k}$ has $k$
components with probability 1. Finally, the time-reversal $(\Pi_n
, \Pi_{n-1}, \ldots , \Pi_1)$ of any fragmentation chain
$(\Pi_1,\ldots,\Pi_n)$, is called a {\em coalescent}.

\medskip A basic problem, only partially solved in this paper, is
the following:

\begin{problem} \label{pbm1}
For which weight sequences $(w_1, \ldots, w_{n-1})$ does there
exist a $\Pnn$-valued Gibbs fragmentation process with these
weights?
\end{problem}

The above definitions were made in terms of $\Pnn$-valued
processes, as this formalism seems most convenient for
computations with Gibbs distributions. Parallel definitions can be
made in terms of $\Pn$-valued processes, using the partial
ordering of refinement on $\Pn$ defined as follows: for partitions
$\la$ and $\mu$ of $n$, $\la$ is a refinement of $\mu$ if and only
if there exist corresponding partitions $\la'$ and $\mu'$ of $[n]$
such that $\la'$ is a refinement of $\mu'$. Less formally, some
parts of $\la$ can be coagulated to form the parts of $\mu$. The
notions of Gibbs distributions and refining sequences transfer
between $\Pn$ and $\Pnn$ in such a way that the following results
can be formulated with either state space. The many-to-one
correspondence between partitions of the set $[n]$ and partitions
of the integer $n$, quantified by (\ref{nfactor}), provides a
many-to-one correspondence between $\Pn$-valued and $\Pnn$-valued
processes, in such a way that the partial ordering of refinement
is preserved. Thus for a $\Pnn$-valued Gibbs fragmentation process
there is a corresponding $\Pn$-valued fragmentation process, and
vice-versa.

\medskip To see that Problem \ref{pbm1} is of some interest, note that
Example \ref{exm3} (cutting a random rooted segment) provides a
$\Pnn$-valued Gibbs fragmentation process with weights $w_j = j!$
for each $n$. Example \ref{exm4} (cutting a random rooted tree)
does the same thing for the weights $w_j = j^{j-1}$. What about
for the sequence $w_j = 1$ of Example \ref{exm1} (uniform random
partitions) or the sequence $w_j = (j-1)!$ of Example \ref{exm2}
(uniform random permutations)? In these examples it is not obvious
how to construct a Gibbs fragmentation process for $n \ge 4$.
(Note that for $n \le 3$ there exists a $\Pnn$-valued Gibbs
fragmentation process for arbitrary positive weights $w_1$ and
$w_2$, for trivial reasons.) The question for $w_j = 1$ is largely
settled by the following proposition, which can be traced as far
back as \cite{harper} (see also \cite{janson}). See Section
\ref{continuous} regarding $w_j = (j-1)!$.

\begin{proposition}\label{nonexist}
There is an $n_0 < \infty$ such that for all $n \ge n_0$ there
does not exist a $\Pnn$-valued Gibbs fragmentation process
$(\Pi_{k}, k \in [n])$ with equal weights $w_1 = w_2 = \cdots =
w_{n-1}$.
\end{proposition}

\proof Let $\Pi_{[n,k]}$ denote a random partition with the Gibbs
distribution on $\Pnk$ with equal weights $w_1 = w_2 = \cdots =
w_{n-1}$, meaning that $\Pi_{[n,k]}$ has uniform distribution on
$\Pnk$. Let
\begin{equation}
\label{rankedxs} X_{(n,k,1)} \ge X_{(n,k,2)} \ge \cdots \ge
X_{(n,k,k)}
\end{equation} denote the
sizes of components of $\Pi_{[n,k]}$ arranged in decreasing order.
Then for each fixed $i$ and $k$ with $1 \le i \le k$ the $i$th
largest component of $\Pi_{[n,k]}$ has relative size
$X_{(n,k,i)}/n$ which converges in probability to $1/k$ as $n \te
\infty$. This follows easily from the law of large numbers, and
the elementary fact that $\Pi_{[n,k]}$  has the same distribution
as $\Pi_{n,k}^*$  given that $\Pi_{n,k}^*$  has $k$ components,
where $\Pi_{n,k}^*$ is the random partition of $[n]$ generated by
$n$ independent random variables $U_1, \ldots, U_n$ each with
uniform distribution on $[k]$. (So $i$ and $j$ are in the same
component of $\Pi_{n,k}^*$ if and only if $U_i = U_j$.) In
particular, there is an $n_0 < \infty$ such that for all $n \ge
n_0$ both
\begin{equation}
\label{5121} \PR( X_{(n,2,2)}
>  (5/12) n ) > 1/2
\end{equation}
and also
\begin{equation}
\label{5122} \PR( X_{(n,3,1)}
> (5/12) n  ) <  1/2
\end{equation}
But if $(\Pi_{[n,k]}, 1 \le k \le n)$ were a fragmentation
process, then $\Pi_{[n,3]} $ would be derived from $\Pi_{[n,2]} $
by splitting one of the two components of $\Pi_{[n,2]}$, and hence
$X_{(n,2,2)} \le X_{(n,3,1)}$ with probability one.
Thus for a fragmentation process, (\ref{5121}) implies the reverse
of the inequality (\ref{5122}), and this contradiction yields the
result.
\endpf

The above argument proves the non-existence for large $n$ of a
Gibbs fragmentation process for any weight sequence $(w_j)$ such
that for $k=2$ or $k=3$ the components in the Gibbs partition of
$[n]$ into $k$ components are approximately equal in size with
high probability. For the weight sequence $w_j = j!$ of Example
\ref{exm3}, what happens instead is that the sequence of ranked
sizes (\ref{rankedxs}), normalized by $n$, has a non-degenerate
limit distribution for each $k$. As observed in \cite[\S
5]{ki82co}, this limit distribution on $[0,1]^k$ is the
distribution of the ranked lengths of $k$ subintervals of $[0,1]$
obtained by cutting $[0,1]$ at $k-1$ points picked independently
and uniformly at random from $[0,1]$. This asymptotic distribution
has been extensively studied \cite{holst80stick}. For the weight
sequence $w_j = j^{j-1}$ of Example \ref{exm4}, the behavior is
different again. What happens is that for each fixed $k$ the
sequence of ranked sizes (\ref{rankedxs}), when normalized by $n$,
converges in probability to $(1,0, \ldots 0)$. That is to say, for
any fixed $k$, for sufficiently large $n$, after $k$ steps in the
fragmentation process, there is with high probability one big
component of relative mass nearly 1, and $k-1$ small components
with combined relative mass nearly zero. To be more precise, it is
easily shown that the $k-1$ small components, when kept in the
order they are broken off the big component, have unnormalized
sizes $X_{n,1}, \ldots, X_{n,k-1}$ that are approximately
independent for large $n$ with asymptotic distribution
\begin{equation} \label{borel distr}
\lim_{n \te \infty } \PR(X_{n,i} = j ) = { j^{j-1} \over j! }
e^{-j} ~~~(i,j = 1,2, \ldots )
\end{equation}
which is the {\em Borel distribution} of the total progeny of a
critical Poisson-Galton-Watson process with Poisson(1) offspring
distribution started with one individual, which can be read from
\re{tprogkn} and \re{poisb}. See \cite[\S 4.1]{randomforests}
for proofs and various generalizations. As a consequence of
(\ref{borel distr}) and the asymptotic independence of the
$X_{n,1}, \ldots, X_{n,k-1}$, the asymptotic distribution of the
combined size $X_{n,1} + \ldots + X_{n,k-1}$ of all but the
largest component of the partition of $[n]$ into $k$ components is
the distribution of the total progeny of the Poisson-Galton-Watson
process with Poisson(1) offspring distribution started with $k$
individuals, which is the {\em Borel-Tanner} distribution
\cite{consul-famoye}
\begin{equation} \lim_{n \te \infty } \PR(X_{n,1} + \ldots +
X_{n,k-1} = m ) = \frac{k-1}m\frac{m^{m-k+1}e^{-m}}{(m-k+1)!}
\end{equation}
which can also be read from \re{tprogkn} and \re{poisb}. According
to the classification of Barbour and Granovsky \cite{bg}, the
examples corresponding to $w_j=j!$, $w_j=j^{j-1}$ and $w_j=(j-1)!$
belong respectively to the expansive, convergent and logarithmic
structures. These structures exhibit quite a different asymptotic
behavior, which may account for the differences observed here
between these three examples.

\medskip
As a contrast to Proposition \ref{nonexist}, it is known
\cite{whittle3,durrett-granovsky-gueron} that for any strictly
positive sequence of weights $(w_j)$, there is a reversible
coagulation-fragmentation process on $\Pnn$ with the canonical
Gibbs distribution (\ref{cangibbs})
as its equilibrium distribution.

\section{Existence of Gibbs fragmentation processes}
\label{main result}

The problem of the existence of a Gibbs fragmentation $(\Pi_1,
\ldots , \Pi_n)$ for a given integer $n$ and weight sequence
$(w_1,\ldots,w_{n-1})$ is one of existence of an increasing
process on a partially ordered set with contraints on the marginal
distributions of the process. In principle, this is solved by the
work of Strassen on measures with given marginals. See for
instance \cite[Theorem 1 and Proposition 4]{Kamae}. According to
this result,
%
for the existence of a Gibbs fragmentation process it is both
necessary and sufficient that for all $A \subset \Pnk$,
\begin{equation} \label{maxflow}
\sum_{\pi \in A} p_{n,k} (\pi)\le \sum_{\pi'\in A'} p_{n,k+1}
(\pi')
\end{equation}
where $A'$ is the set of partitions $\pi'$ that can be obtained by
splitting a single block of some partition $\pi \in A$, and
$p_{n,k}$ is the Gibbs measure on partitions of $[n]$ into $k$
blocks. A variation of this condition can also be given in terms
of integer partitions rather than set partitions. Unfortunately,
it seems hard to use this general criterion to prove any existence
result. However (\ref{maxflow}) can be used as an algorithm to
check the existence of a Gibbs fragmentation for a finite $n$ and
a given sequence $(w_1,\ldots,w_{n-1})$.
For instance, when $w_j=1$ it is possible to check that the first
$n$ for which the existence of a Gibbs fragmentation fails is
$n=20$, as mentioned in \cite{harper}.

If there exists any $\Pnn$-valued Gibbs fragmentation process
governed by $(w_1, \ldots, w_{n-1})$, then there exists one that
is a Markov chain. For given a non-Markovian process, one can
always create a Markov chain with the same one-step transition
probabilities and the same marginal distributions. So Problem
\ref{pbm1} reduces to:

\begin{pbm} {pbm2}
For which weight sequences $(w_1, \ldots, w_{n-1})$ does there
exist a transition matrix $\{P(\pi , \nu )\}$ indexed by $\Pnn$
such that $P(\pi , \nu ) > 0 $ only if $\nu$ is a refinement of
$\pi$, and if $\nu \in \mathcal{P}_{[n,k]}$
\begin{equation} \sum_{\pi \in \Pnn}
p_{n,k-1}(\pi) P( \pi, \nu) = p _{n,k} (\nu  ) ~~~(2 \le k \le n)
\end{equation}
where $p_{n,k}( \nu )$ is given by the microcanonical Gibbs
formula {\em (\ref{gibbsnk})} ?
\end{pbm}

Such a transition matrix $P(\pi, \nu)$ corresponds to a {\em
splitting rule} which describes for each $1 \le k \le n-1$ and
each partition $\pi$ of $[n]$ into $k-1$ components, the
probability that $\pi$ splits into a partition $\nu$ of $[n]$ into
$k$ components.  Given that $\Pi_{k-1} = \pi_{k-1}$ with
$\pi_{k-1} = \{A_1', \ldots, A_{k-1}'\}$ say, the only possible
values $\pi_k$ of $\Pi_k$ are those $\pi_{k} = \{A_1, \ldots,
A_{k}\}$ such that two of the $A_j$
form a partition of one of the $A_i'$, and the remaining $A_j$ are
identical to the remaining $A_i'$. The initial splitting rule
starting with $\pi_1 = \{[n]\}$ is assumed to be specified by the
Gibbs distribution $p_{n,2}$ determined by the weight sequence
$(w_1, \ldots, w_{n-1})$ for $n_1$ and $n_2$ with $n_1 + n_2 = n$.
That is from (\ref{gibbsnk}):
\begin{equation}
\label{gibbs2} \PR( \Pi_2 = \{A_1, A_2 \} )
= \frac{ w _{\# A_1} w_{\# A_2 } } { B_{n,2}( w_1, \ldots,
w_{n-1}) }
\end{equation}
for $B_{n,2}$ as in (\ref{belldef}). The simplest way to continue
is to use the following:

\mn {\bf Recursive Gibbs Rule:} for each $1 \le k \le n-1$, given
that $\Pi_{k-1}= \{A_1', \ldots, A_{k-1}'\}$ and that some
particular block $A \in \{A_1', \ldots, A_{k-1}'\}$ is split with
$\# A = m \ge 2$, $\Pi_k$ is obtained by splitting this block is
split into $\{ A_1 , A_2 \}$ with probability given by the right
side of \re{gibbs2} with $m$ instead of $n$.

To complete the description of a splitting rule, it is also
necessary to specify for each partition $\pi_{k-1} = \{A_1',
\ldots, A_{k-1}'\}$ the probability that the next component to be
split is $A_i'$, for each $1 \le i \le k-1$. The simplest possible
assumption seems to be the following:

\mn {\bf Linear Selection Rule:} Given $\pi_{k-1} = \{A_1',
\ldots, A_{k-1}'\}$, split $A_i'$ with probability proportional to
$\# A'_i - 1$, that is with probability $(\# A'_i - 1)/(n-k+1)$.

\medskip While this selection rule is somewhat arbitrary, it is
natural to investigate its implications for the following reasons.
Firstly, blocks of size $1$ cannot be split, so the probability of
picking a block to split must depend on size. This probability
must be 0 for a block of size 1, and 1 for a block of size
$n-k+1$. The simplest way to achieve this is by linear
interpolation. Secondly, both the segment splitting model and the
tree splitting model described in Examples \ref{exm3} and
\ref{exm4} follow this rule. In each of these examples a block of
size $m$ is a graph component with $m-1$ edges, so the linear
selection rule corresponds to picking an edge uniformly at random
from the set of all edges in the random graph whose components
define $\Pi_{k-1}$. Given two natural combinatorial examples with
the same selection rule, it is natural to ask what other models
might follow the same rule.

More complex splitting rules are also of interest. Consider for
instance a continuous time Markov fragmentation chain which
fragments blocks of size $j$ according to infinitesimal rates
dictated by $\lambda_j p_{j,2}$ for some $\lambda_j >0$. The
embedded discrete time chain is then a recursive Gibbs
fragmentation chain with the property that the probability to
select a particular block of size $n_i$ for the next fragmentation
is proportional to $\lambda_{n_i}$.  But we do not know any nice
description of the law of $\Pi_k$ in this case beyond saying that
it is the solution of some Kolmogorov forward equations. See also
\cite{hmpw} for theory of discrete and continuous time Markov
fragmentation chains which are exchangeable and consistent as $n$
varies, meaning that they can be associated with fragmentations of
a mass continuum.

\subsection{Main results}

This section presents the main results, whose proofs are provided
in the next two sections. Recall the definition of the discrete
Marcus-Lushnikov coalescent process on $\Pnn$ with affine kernel:
this is the unique Markov chain on $\Pnn$ such that $\pi_1$ is the
partition consisting of singletons and $\pi_k$ is obtained from
$\pi_{k-1}$ by merging each pair of blocks of sizes $i$ and $j$
with probability proportional to $K_{i,j} =a+b(i+j)$ for some
constants $a$ and $b$. In the case $a=1$ and $b=0$ this is
Kingman's $n$-coalescent, as described in Example 3 (blocks
coalesce at rate 1), while if $a=0$ and $b=1$ this is the additive
coalescent mentioned in Example 4.

\begin{theorem}\label{new}
Fix $n \ge 4$, and let $(w_j,1 \le j \le n-1)$ be a sequence of
positive weights with $w_1 = 1$. The following two statements are
equivalent:

\begin{enumerate}  \item[{\em (i)}] \nopagebreak
The $\Pnn$-valued fragmentation process $(\Pi_k, 1 \le k \le n)$
defined by the recursive Gibbs splitting rule derived from these
weights, with the linear selection rule, is such that for each $1
\le k < n$ the random partition $\Pi_k$ has the microcanonical
Gibbs distribution $p_{n,k}$ with the same weights.

\item[{\em (ii)}] The weight sequence $w_j$ is of the form
\begin{equation} \label{wjbc}
 w_j = w_j^{b,c} :=
\prod_{i = 2}^j ( i c  + j b )
\ , \ (j=2,\ldots, n-1)
\end{equation}
for some real $b$ and $c$ such that \eq \label{bcineq} \mbox{$b +
c >0$ and either $b\ge 0$ or $b < 0$ and $c > -(n-1)b/2$.} \en
\end{enumerate}
\begin{enumerate}
\item[{\em (iii)}] The time reversal of $(\Pi_k, 2 \le k\le n)$ is
a discrete Marcus-Lushnikov coalescent with affine kernel
$K_{i,j}=a+b(i+j)$. Moreover, in this case the $b$ is the same as
in {\em (ii)} and $a=2c$.
\end{enumerate}
\end{theorem}

Note that $c$ and $b$ appearing in (ii) and (iii) are unique only
up to a constant common factor.

Given a continuous time $\Pnn$-valued coalescent or fragmentation
process $(\Pi(t), t \in I)$, define the {\em discrete skeleton} of
$(\Pi(t), t \in I)$ to be the $\Pnn$-valued process $(\Pisk{k}, 1
\le k \le n)$ where $\Pisk{k}$ is the common value of $\Pi(t)$ for
all $t \in I$ such that $\#\Pi(t) = k$. Provided either $b\ge 0$
or $b < 0$ and $c > -nb/2$ the time-reversed process in part (iii)
of the above theorem is the discrete time skeleton of the
continuous time affine coalescent with collision rate kernel
$K_{i,j}:= 2 c + b( i + j )$. As observed by Hendriks et al.
\cite{hses85}, this kernel has the special property that the
process $(\# \Pi(t), t \ge 0 )$ is independent of the discrete
skeleton of $(\Pi(t),t\ge 0)$ (in the special case of Kingman's
coalescent this had been proved earlier in \cite{ki82co}). Thus
Theorem \ref{new} implies the result of Hendriks et al.
\cite{hses85} that for an affine coalescent in continuous time the
distribution of $\Pi_t$ is a Gibbs distribution with weights
$w_j^{b,c}$ as in (\ref{wjbc}), that is a mixture over $k$, with
mixing weights depending on $t$, of the microcanonical Gibbs
distributions $p_{n,k}^{b,c}$ featured in Theorem \ref{new}. The
fact that $(\Pi(t), t \ge 0)$ is a Gibbs coalescent with a
particular sequence of weights $w_j$ is related in this instance
to the fact that its discrete skeleton is a Gibbs coalescent with
the same weights.   But this equivalence relies on the
independence of the process $(\# \Pi(t), t \ge 0 )$ and its
discrete skeleton. It is not always true that the discrete
skeleton of a continuous time Gibbs coalescent is a discrete time
Gibbs coalescent, as illustrated by example in Section
\ref{continuous}.

The following corollary was suggested by comparison of Theorem
\ref{new} with the branching process interpretation of Bell
polynomials provided in Section \ref{branching}. To simplify the
argument we introduce a regularity condition on the offspring
distribution \re{simplify}, but this assumption may not be
strictly necessary for the result to stay valid.

\begin{corollary}
\label{forests} Fix $n \ge 4$. Let $T$ denote a Galton-Watson tree
with offspring distribution $(p_j)$ such that \eq \label{simplify}
\mbox{ $p_j >0$ if and only if $0 \le j < j_1$ for some $1 \le j_1
\le \infty$. } \en Let $F_1$ be $T$ conditioned to have $n$ nodes,
regarded as a random plane tree (a tree with ordered branches),
and for $2 \le k \le n$ let $F_k$ be the plane forest of $k$ trees
obtained by first cutting $k-1$ edges of $F_1$ picked by a process
of random sampling without replacement, and then putting these $k$
trees in random order, with all $k!$ orders equally likely, where
the cutting and ordering processes are independent of each other
and of $F_1$. Then the following conditions are equivalent:
\begin{itemize}
\item[{\em (i)}] The offspring distribution is such that \eq
\label{off} \frac{ p_j }{p_0}  = \frac{1}{j!} \prod_{i = 1}^j ( b
- ( i - 2 ) c  ) \en for some real parameters $b$ and $c$ with
$b+c >0$ and such that the product is non-negative for all $1 \le
j \le n-1$.

\item[{\em (ii)}] The forest of two trees $F_2$ is distributed
like two independent copies of $T$ conditioned to have a total of
$n$ nodes.

\item[{\em (iii)}] For every $1 \le k \le n$ the forest of $k$
trees $F_k$ is distributed like $k$ independent copies of $T$
conditioned to have a total of $n$ nodes.

\end{itemize}

For such a sequence of forests $(F_k, 1 \le k \le n)$ let $(\Pi_k,
1 \le k \le n)$  be the refining sequence of partitions of $[n]$
defined by labelling the $n$ nodes of tree $F_1$ by a random
permutation independent of $F_1$, and letting the blocks of
$\Pi_k$ be the tree components of $F_k$. Then

\begin{itemize}
\item[{\em (iv)}] the sequence of partitions $(\Pi_1, \ldots,
\Pi_n)$ develops by recursive Gibbs fragmentation  with linear
selection, for the weight sequence $(w_j^{b,c})$ as in {\em
(\ref{wjbc})}, and $(\Pi_n, \ldots, \Pi_1)$ is a Marcus-Lushnikov
coalescent with the affine kernel $K_{i,j}=2c+b(i+j)$.
\end{itemize}
\end{corollary}

The implications (i) $\Rightarrow$ (ii) $\Rightarrow$ (iii)
$\Rightarrow$ (iv) of this Corollary were provided in
\cite{randomforests} for the case $c=0$, when the offspring
distribution can be Poisson with mean $b$ for any $b >0$. Note
that (i) only specifies the conditional offspring distribution
given at most $n-1$ children, as is necessary for the converse for
a fixed $n$. The conditions on $b$ and $c$ imposed in (i), which
are necessary for construction of the forest-valued fragmentation
$(F_k)$, imply but are not implied by the conditions \re{bcineq}
which are necessary for construction of the partition-valued
fragmentation $(\Pi_k)$. To illustrate for $n =4$, the conditions
\re{bcineq} are that $b + c >0$ and $3b + 2 c > 0$, whereas those
in (i) above are $b + c >0$ and either $b - c \ge 0$ or $b = 0$.
In either case, $b \ge 0$, hence $3b + 2 c >0$, but not
conversely. The $b$ and $c$ such that the conditions \re{bcineq}
hold for all $n$ are those with $b \ge 0 $ and $b + c >0$. Whereas
there is the forest-valued representation for all $n$  if and only
if one of the following further conditions holds, as discussed
later in Section \ref{branching}.
\begin{itemize}
\item $c = 0$: the offspring distribution is then Poisson$(b)$;
\item $c >0$ and $b = (a-1)c$ for a positive integer $a$: the
offspring distribution is then binomial$(a,p)$ for $p = c/(c+1)$;
\item $-1< c < 0$ and $b= (a-1)c$ with $-a = r>0$: the offspring
distribution is then negative binomial$(r,p)$ for $p = c+1$.
\end{itemize}
The cases when $a$ is an integer admit further combinatorial
interpretations, which we will discuss in more detail elsewhere.
For instance, when $a=-1$ we obtain a representation of the affine
coalescent with collision kernel $K_{i,j} = i+j -1$ by time
reversal of a process of coalescent plane forests $(F_k, 1 \le k
\le n)$, where the forest with $k$ trees has the uniform
distribution on the set of
$$
\frac{k}{n} { 2 n - k - 1 \choose n - k }
$$
plane forests with $k$ trees. And when $a$ is a positive integer,
there is an interpretation of the $(a - 1)(i+j) + 2$ coalescent in
terms of trees where each node has either $0$ or $a$ children.

\section{Preliminaries}

\subsection{Generating functions}
\newcommand{\winv}{w^{(-1)}}
\newcommand{\hinv}{h^{-1}}
\newcommand{\CF}{\widehat{C}}
\newcommand{\FF}{C}

Let $w(z):= \sum_{n=1}^\infty w_nz ^n/n!$ be the exponential
generating function associated with the sequence of weights
$(w_n)$. It follows easily from \re{belldef} that \begin{equation}
\label{bnk1} B_{n,k}(w_1,w_2, \ldots) = \frac{n!}{k!} [z^n] w(z)^k
\end{equation} where $[z^n] w(z)^k$ denotes the coefficient of
$z^n$ in the expansion of $w(z)^k$ in powers of $z$. In
particular, \begin{equation} \label{bn2} B_{n,2}(w_1, \ldots,
w_{n-1}) = \frac{1 }{ 2 } \sum_{l=1}^{n-1} {n \choose l } w_l
w_{n-l} . \end{equation} Assuming the weights are such that
$w(\xi) < \infty$ for some $\xi >0$, the formula $$
\PR(Y=n)=\frac{w_n \xi^n}{n!w(\xi)} $$ defines the distribution of
a non-negative random variable $Y$ whose probability generating
function is \begin{equation} \label{kolgf} E ( z^Y ) = w(  \xi z
)/w(\xi) . \end{equation} If $Y_1,Y_2, \ldots$ is a sequence of
independent random variables with the same distribution as $Y$,
then \begin{equation} \label{kolgf1} \PR( Y_1 + \cdots Y_k = n ) =
[z^n] \left( \frac{ w ( z \xi ) }{ w (\xi } \right)^k = \frac{ k!
B_{n,k} \xi ^n }{ n! w(\xi)^k }, \end{equation} which appears for
instance in (1.3.1) of \cite{kolchin rg} and \cite[Lemma 3.1]{eg}.
This implies the {\em Kolchin representation of block sizes in a
Gibbs partition} \cite[Theorem 1.2]{saintflour}: for a random
partition of $[n]$ with the microcanonical Gibbs distribution
$p_{n,k}$ derived from $(w_j)$, when the $k$ blocks are put in a
random order, with each of $k!$ possible orders equally likely,
independently of the sizes of the blocks, the sequence of block
sizes is distributed as \eq \label{kolchin} (Y_1,\ldots, Y_k)
\text{ given }  Y_1 + \cdots Y_k = n . \en Easily from \re{bnk1}
there is the {\em exponential formula} \eq \label{expform} e^{x
w(z)} = \sum_{n = 0}^\infty C_n(x) z^n \en where $C_0(x) = 1$ and
$C_n(x)$ for $n = 1,2, \ldots$ is the polynomial $$ C_n(x) =
(n!)^{-1} \sum_{k = 1}^n B_{n,k}(w_1,w_2, \ldots) x^k . $$ The
polynomials  $C_n(x)$ are then of {\em convolution type}, meaning
that for $n\ge 1$, \eq \label{convtype} C_n(x+y) = \sum_{k = 0}^n
C_{k}(x) C_{n-k}(y) . \en Assuming now that $w_1 = 1$, let $\winv$
denote the compositional inverse of $w$ defined by $\winv (w(z)) =
z$. According to the Lagrange inversion formula \cite[Theorem
5.4.2]{stanleyv2} \begin{equation} \label{lagrange} [z ^n ] w(z)^
k = \frac{k}{n} [ z^{n-k}] \left( \frac{z } { \winv(z) } \right)^n
= \frac{k}{n} \CF_{n-k} (n) \end{equation} where $\CF_n(x)$ is the
sequence of polynomials of convolution type defined by
\begin{equation} \label{bnk22} \left( \frac{z} { \winv(z) }
\right)^x = \sum_{n=0}^\infty \CF_n(x) z^n. \end{equation}
Combining
 \re{bnk1} and \re{lagrange} we obtain the following lemma. See also
Knuth \cite{knuth} for a similar discussion.

\begin{lemma}
\label{bellrep} Each sequence of real weights $(w_1, w_2, \ldots)$
with $w_1 = 1$ admits the representation
\begin{equation}
\label{bnk33} w_n = (n-1)! \, \CF_{n-1}(n)
\end{equation}
for a unique sequence of polynomials $\CF_n(x)$ of convolution
type, namely that determined by {\em \re{bnk22}}, in which case
for $n\ge 1$,
\begin{equation}
\label{bnk34} B_{n,k}(1, w_2, w_3, \ldots ) = \frac{ (n-1)! }{ ( k
- 1)! } \CF_{n-k}(n) .
\end{equation}
\end{lemma}

Many sequences of polynomials of convolution type are known
\cite[Examples 2.2.16]{dibucch}, each providing a sequence of
weights $(w_n)$ for which the Bell polynomials can be explicitly
evaluated using \re{bnk34}. As a general rule, weight sequences
with manageable formulas for the $B_{n,k}$ are those with a simple
formula for $z/\winv(z)$ rather than for $w(z)$.  A rich source of
such examples is provided by the theory of Galton-Watson branching
processes.

\subsection{Branching Processes}
\label{branching} Given a weight sequence $(w_j)$ with $w_1=1$ and
exponential generating function $w(z) = \sum_{n\ge 1} w_n z^n/n!$,
let \eq G( z ):= \frac{ z }{ \winv(z) } = \sum_{n = 0}^\infty
\CF_n(1) z^n . \en where $\CF_0(1) = 1$. Then provided \eq
\CF_n(1) \ge 0 \mbox { for all $n\ge 1$ and } G(\eta) < \infty
\mbox{ for some } \eta
> 0 \en the formula \eq \label{gG} g(z):= \frac{G(z
\eta)}{G(\eta)} = \sum_{n = 0}^\infty \frac{ \CF_n(1) \eta
^n}{G(\eta)} \, z^n \en defines the probability generating
function of a non-negative integer valued random variable $X$ with
distribution \eq \label{offspring} \PR(X = n ) = \frac{ \CF_n(1)
\eta ^n}{G(\eta)}  ~~~~~(n = 0,1,2, \ldots ) \en Conversely, for
each distribution of $X$ with $\PR(X = 0 ) >0$ and each  $\eta >0$
it is easily seen that there is a unique sequence of convolution
polynomials $\CF_n$ such that \re{offspring} holds.  This is a
particular case of \cite[Theorem 2.1.14]{dibucch}. Let $Y$ be the
total progeny in a Galton-Watson branching process with generic
offspring variable $X$. It is well known that the probability
generating function \eq \label{totalprog} h(z):=  \sum_{n =
1}^\infty \PR(Y=n) z^n \en can be characterized as the unique
solution of the functional equation \eq \label{hg} h(z) = z
g(h(z)) \en which is obtained by conditioning on the number of
offspring of the root individual \cite[Section 6.1]{saintflour}.
Note in particular that given the generating function $h$ of the
total progeny, the offspring probability generating is determined
by
\begin{equation}
\label{g given h} g(v)=\frac{v}{h^{(-1)}(v)}
\end{equation}
Also,
$$
\label{proper} \mbox{ $h(1) = 1$, meaning $\PR(Y < \infty) = 1$, }
$$
if and only if mean of the offspring distribution is at most $1$,
that is by (\ref{gG}) \eq \label{mucond} g'(1) = \eta
G'(\eta)/G(\eta) \le 1. \en Note that the assumed form
\re{offspring} of the offspring distribution forces $\PR(X=0) >0$
(since $\hat C_0(1)=1$), and so forbids the degenerate case with
$\PR(X=1) = 1$. Combining this discussion with Lemma \ref{bellrep}
we obtain:

\begin{proposition}
\label{bprep} Let $w_1 = 1, w_2, \ldots$ be a sequence of
non-negative weights with exponential generating function $w(z):=
\sum_{n = 1}^\infty w_n z^n /n!$ and let $\winv$ be the
compositional inverse of $w$ defined by $w^{-1}(w(z))=z$. The
following two conditions are equivalent:

\begin{enumerate}
\item[{\em (i)}] there exists $\xi >0 $ such that $w(\xi) <
\infty$ and the random variables $Y_i$ in Kolchin's representation
{\em \re{kolchin}} of Gibbs partitions, with generating function
$w(z\xi)/w(\xi)$, are distributed like the total progeny of some
Galton-Watson branching process started with one individual.

\item[{\em (ii)}] The power series
$$
G(z):= \frac{z}{\winv(z)} = \sum_{n = 0}^\infty \CF_n(1) z^n
$$
has non-negative coefficients $\CF_n(1)$ and $G(\eta) < \infty$
for some $\eta >0$.
\end{enumerate}
When these conditions hold, the offspring distribution is as
displayed in {\em (\ref{offspring})}, with generating function
$g(z) = G( \eta z )/G(\eta)$ for $\eta = w(\xi)$, and $g$ must
satisfy $g'(1) \le 1$. The associated evaluation of Bell
polynomials is  then \eq \label{belleval} B_{n,k} (1, w_2, w_3,
\ldots ) = \frac{ (n-1)! }{ ( k - 1)! } \CF_{n-k}(n) = \frac{ n!
w(\xi)^k }{k! \xi^n} \PR( Y_1 + \cdots +Y_k = n ) \en where
$\CF_n(x):= [z^n] G(z)^x$ and $Y_1 + \cdots +Y_k$ represents the
total progeny of the branching process started with $k$
individuals.
\end{proposition}
\proof Condition (i) is that $w(z \xi)/w(\xi) = h(z)$ where $h$ is
derived from some probability generating function $g$ via \re{hg}.
Let $\hinv$ denote the compositional inverse of $h$, defined by
$h(\hinv(z)) = z$. The equation $h(z) = v$ is solved by $z =
\winv( v w(\xi ) ) /\xi$, so
using (\ref{g given h}) $g$ is recovered as
$$
g( v ) = \frac{v}{\hinv(v)} = \frac{ v \xi } { \winv( v w (\xi ) )
} = \frac{ \xi } { w(\xi ) } \frac { v w(\xi ) }{ \winv ( v w(\xi
) ) } = \frac{ G( \eta v ) }{ G ( \eta ) }
$$
where $\eta = w(\xi)$ so that $G(\eta) = G( w ( \xi ) ) = \frac{
w( \xi  ) } { \winv ( w ( \xi ) )} = \frac{ w ( \xi ) }{\xi }.$
The rest is read from Lemma \ref{bellrep}.
\endpf

The conditions of the previous proposition force the branching
process to be critical or subcritical. For arbitrary $\eta$ with
$G(\eta) < \infty$, and a branching process with offspring
generating function $g(z):= G(z \eta)/G(\eta)$, the Lagrange
inversion formula shows that the distribution of the total progeny
of the branching process started with $k$ individuals is given by
the formula
\begin{equation}
\label{tprogkn} \PR( Y_1 + \cdots + Y_k = n ) = [z]^n h(z)^k =
\frac{k}{n} [ z^{n - k } ] g (z ) ^n = \frac{ \eta^{n-k} }{ G(\eta
) ^n } \frac{k}{n} \CF_{n-k}(n)
\end{equation}
where the $Y_i$ are independent and identically distributed
according to this formula for $k=1$. Formula \re{tprogkn} can be
rewritten using Lemma \ref{bellrep} as
\begin{equation}
\label{tprogbell} \PR( Y_1 + \cdots Y_k = n ) = \frac{ \xi^n  }{
\eta^{k} } \frac{ k ! }{n!} B_{n,k} (1,w_2, w_3, \ldots)
\end{equation}
which is also consistent with (\ref{kolgf1}). Here $\xi:=
\eta/G(\eta)$, and necessarily $w(\xi) \le \eta$, with $w(\xi) =
\eta$ and $\PR( Y_1 + \cdots + Y_k < \infty) = 1$ only in the
critical or subcritical case $g'(1) \le 1$.

To illustrate these results, consider first the generating
function $G(z)= e^{b z}$ so that
$$
G(z)^x = e^{b z x } = \sum_{n = 0}^\infty { b^n  x^n \over n! }
z^n
$$
The associated sequence of convolution polynomials is
$$
\CF_n(x) = b^n x^n /n! .
$$
The convolution identity \re{convtype} is the binomial theorem.
The corresponding weight sequence is
$$
w_n = (n-1) ! \CF_{n-1}(n) = b^{n-1} n^{n-1}
$$
and the Bell polynomial evaluation is \eq \label{poisb} B_{n,k} =
\frac{(n-1)!}{(k-1)!} \CF_{n-k}(n) = {n - 1 \choose k-1 } b^{n-k}
n^{n-k} \en as indicated earlier in \re{cayley}. The branching
process interpretation is that for Poisson offspring distribution
with mean $b$, the distribution of the total progeny of the
branching process started with $k$ individuals is given by formula
\re{tprogkn} with the above substitutions for $\eta = 1$ and $\xi
= 1/G(1) = e^{-b}$.

Consider next the generating function $G(z)= ( 1 + c z )^a$ for
some pair of real parameters $a$ and $c$, so that
$$
G(z)^x=  (1 + c z )^{ax} = \sum_{n = 0}^\infty { a x \choose n }
c^n z^n.
$$
The associated sequence of convolution polynomials is
$$
\CF_n(x) = { a x \choose n } c^n  .
$$
In this case, the convolution identity \re{convtype} is
called the Chu-Vandermonde identity (see, e.g., \cite{gould}). The
corresponding weight sequence is \eq \label{nbwts} w_n = (n-1)
!\CF_{n-1}(n) = (n-1)! {a n \choose n - 1} c^{n-1} \en and the
Bell polynomial evaluation is \eq \label{bpac} B_{n,k} =
\frac{(n-1)!}{(k-1)!} \CF_{n-k}(n) = \frac{(n-1)!}{(k-1)!} { a n
\choose n - k } c^{n-k} . \en Two cases of this formula have well
known probabilistic interpretations
\cite{consul-shenton,consul-famoye,devroye}, as indicated in the
next two paragraphs.

If $a$ is a positive integer and $c > 0 $,
then $\CF_n(1) \ge 0$ for all $n$. For $\eta = 1$ the probability
generating function \re{gG} is
$$
g(z) = \frac{G(z)}{G(1)}  = \left( \frac{ 1 + c z }{ 1 + c }
\right)^a = (q + pz)^a
$$
for $p:= c/(1+c)$ and $q:= 1-p$. This represents the binomial
distribution with  parameters $a$ and $p$. For $a=1$ the
evaluation \re{bpac} reduces to  the previous evaluation \re{lah}
of the Lah numbers. The branching process in this case is a rather
trivial one, with each individual having either $0$ or $1$
offspring. So the random family tree is just a line of vertices
whose length is geometrically distributed. Cutting the edges in
such a segment of random length by an independent process of
Bernoulli trials yields a geometrically distributed number of
components, which given their number have independent and
identically distributed lengths with another geometric
distribution. According to Corollary \ref{forests} a similar
interpretation of the microcanonical Gibbs distributions with
weights \re{nbwts} can be provided in terms of random cutting of
edges of a Galton-Watson tree both in the case of binomial$(a,p)$
offspring distribution for $a = 1, 2, 3, \ldots$, and in the
following case of negative binomial offspring distribution.

If $a = -r$ and $c = - q$ for $r >0$ and $0 < q < 1$, again
$\CF_n(1) \ge 0$ for all $n$. For $\eta = 1$ the probability
generating function \re{gG} is
$$
g(z) = \frac{G(z)}{G(1)}  = \left( \frac{ 1 + c z }{ 1 + c }
\right)^{-r} = \left( \frac{ 1 - q }{ 1 - q z } \right)^{r}
$$
which is the generating function of the negative binomial
distribution with parameters $r>0$ and $p = 1 - q = 1 + c \in
(0,1)$.

It is easily seen that the coefficients ${ a \choose n } c^n$ are
non-negative for all $n$ only in the two cases just discussed. So
only in these cases does the Bell polynomial \re{bpac} admit the
interpretation  of Lemma \ref{bprep} in terms of the total progeny
of a branching process for all $n$. Still, the weights $w_n$ in
\re{nbwts} are non-negative for other choices of real $a$ and $c$,
for instance $a >1$ and $c >0$. These weights still define a Gibbs
distribution on partitions, and there is the Kolchin
representation \re{kolchin} for the sizes of blocks of such a
partition. A natural probabilistic construction of such random
partitions is provided by Theorem \ref{new}. The interesting
intermediate case, when the coefficients ${ a \choose j } c^j$ are
non-negative only for $ j < j_1$ for some $j_1 < \infty$,
corresponds to Corollary \ref{forests}. Then $\CF_n(1)$ can be set
equal to $0$ for $j \ge j_1$, and the previous branching process
formulas remain valid provided $n$ is restricted to $n \le j_1$.

\subsection{Evaluation of a Bell polynomial}
The results of the last two Bell polynomial evaluations \re{poisb}
and \re{bpac}, which are implicit in the standard theory of
branching processes, are unified algebraically by the following
lemma. The evaluation \re{bellbc} is also implicit in
\cite[(19)-(21)]{hses85}, and plays a key role in our treatment of
Gibbs models for fragmentation processes.

\begin{lemma}
For each pair of real parameters $b$ and $c$, the polynomials \eq
\CF_n^{b,c}(x):= \frac{1 }{n!} \prod_{j = 0}^{n-1} ( b x  + c x -
c j ) \en are of convolution type. For the corresponding weight
sequence
\begin{equation}
\label{bellbc2} w_{n}^{b,c} := (n-1)! \CF_{n-1}^{b,c}(n) =
\prod_{i = 2}^{n} ( i c + n b )      
\end{equation}
there is the Bell polynomial evaluation
\begin{equation}
\label{bellbc} B_{n,k}(1,w_{2}^{b,c},w_{3}^{b,c}, \ldots) = (n-k)!
\CF_{n-k}^{b,c}(n) = {n-1 \choose k - 1 } \prod_{i = k+1}^{n} ( ic
+ n b ). 
\end{equation}
\end{lemma}
\proof This is read from the previous example with generating
function $G(z) = (1 + cz)^a$ for $a = b/c + 1$. The limiting case
$c = 0$ corresponds to $G(z) = e^{zb}$.
\endpf

It is convenient to record here as well an immediate consequence
of \re{bellbc}:

\begin{lemma}
\label{uniquew} The sequence of weights $w_n = w_n^{b,c}$ is the
unique solution of the recursion $w_1 = 1$, $w_2 = 2 b + 2 c$, and
\begin{equation}
\label{wlem} w_n = \frac{2c+nb}{(n-1)} B_{n,2}(w_1, \ldots,
w_{n-1})      ~~~~~~(n  = 2, 3, \ldots ) \en for $B_{n,2}$ as in
{\em \re{bn2}}.
\end{lemma}

\section{Proofs}

\subsection{Proof of the main result}

\label{proof1}

The proof of Theorem \ref{new} is based on the next two lemmas.
\begin{lemma} \label{mainlem} Fix $n \ge 4$ and $3 \le k \le n-1$,
and let $(\Pi_{k-1},\Pi_k)$ be a pair of random partitions of
$[n]$ such that $\Pi_{k-1}$ is distributed according to the
microcanonical Gibbs distribution $p_{n,k-1}$ with weights $w_1
=1,w_2, \ldots, w_n$, and $\Pi_k$ is derived from $\Pi_{k-1}$ by
the recursive Gibbs splitting rule with these weights, and the
linear selection rule. The following two conditions are
equivalent: \begin{enumerate} \item[{\em (i)}] \nopagebreak
$\Pi_k$ has the microcanonical Gibbs distribution $p_{n,k}$ with
the same weights. \item[{\em (ii)}] \nopagebreak The function
\begin{equation} \label{fmdef} f(m):= { (m - 1 ) w_m \over B_{m,2}
(w_1, \ldots, w_{n-1}) } ~~~~~~(2 \le m \le n -1 ) \end{equation}
satisfies \begin{equation} \label{condf} \sum_{1 \le i < j \le k}
\,\,\, f( n_i + n_j ) = g(n , k ) \end{equation} for all sequences
of $k$ positive integers $(n_1, \ldots, n_k)$ with $ \sum_{i=1}^k
n_i = n $ and some function $g(n,k)$. \end{enumerate} When these
conditions hold, \begin{equation} \label{bnrec} g(n,k) = \frac{( n
- k + 1 ) B_{n, k-1} } { B_{n,k}} \end{equation} and the reverse
transition from $\Pi_k$ to $\Pi_{k-1}$ is governed by the
Marcus-Lushnikov coagulation mechanism with kernel $K_{i,j} = f(i
+ j)$. In the case $k = 3$ these conditions are equivalent to \eq
\label{fmbc} f(m) = 2 c + m b \mbox{ for all } 2 \le m \le n - 1
\en and hence to $w_j = w_j^{b,c}$ as in {\em (\ref{wjbc})}, for
some real $b$ and $c$. \end{lemma}

\proof Let $\pi_k$ denote any particular partition of $[n]$ into
$k$ blocks, say $\{A_1, \ldots, A_k\}$ with $\#A_i = n_i, 1 \le i
\le k$. For $1 \le i < j \le n$ let $\pi_{k-1}^{i,j}$ be the
partition of $[n]$ into $k-1$ blocks derived from $\{A_1, \ldots,
A_k\}$ by merging of $A_i$ and $A_j$. The hypothesis of the lemma
implies that
\begin{equation} \label{arg1}
\PR(\Pi_{k-1} = \pi_{k-1}^{i,j}, \Pi_k = \pi_k ) = \PR(\Pi_{k-1} =
\pi_{k-1}^{i,j}) \,{ (n_i + n_j - 1) \over (n - k +1 )} { w_{n_i}
w_{n_j } \over B_{n_i + n_j ,2} }
\end{equation}
and that
\begin{equation} \label{arg2} \PR(\Pi_{k-1} = \pi_{k-1}^{i,j})  = {
w_{n_i + n_j } \prod_{l=1}^k w_{n_l} \over B_{n,k-1} \, w_{n_i }
w_{n_j} }
\end{equation}
Substituting (\ref{arg2}) into (\ref{arg1}) gives
\begin{equation} \label{arg3} \PR(\Pi_{k-1} = \pi_{k-1}^{i,j}, \Pi_k
= \pi_k ) =
\frac{ f( n_i + n_j ) \prod_{l=1}^k w_{n_l} } { (n-k+1 ) B_{n,k-1}
}
\end{equation}
for $f$ derived from the weights as in \re{fmdef}. Summing this
probability over all possible choices of $(i,j)$ with $1 \le i < j
\le k$ yields $\PR(\Pi_k = \pi_k)$, so the equivalence of
conditions (i) and (ii) is clear. Assuming these conditions hold,
\re{bnrec} follows at once: dividing \re{arg3} by the Gibbs
formula for $\PR(\Pi_k = \pi_k)$ gives
$$
\PR(\Pi_{k-1}=\pi_{k-1}^{i,j}|\Pi_k=\pi_k) = f(n_i+n_j)/g(n,k).
$$
as claimed. In the case $k=3$, we deduce \re{fmbc} from the
following lemma, and the weights are then determined by Lemma
\ref{uniquew}.
\endpf

\begin{lemma}
\label{linear f lemma} Fix $n \ge 3$ and let $(f(m), 2 \le m \le
n-1)$ be a sequence such that for every triple of positive
integers $(n_1,n_2,n_3)$ with $n_1 + n_2 + n_3 = n$
\begin{equation} \label{eqnf} f(n_1 + n_2 ) + f(n_2 + n_3 ) +
f(n_1 + n_3 )  = C
\end{equation} for some constant $C$. Then there exist constants
$b$ and $c$ such that $f(m) = 2 c + m b $ for every $2 \le m \le
n-1$, and $C = 2 (3 c + n b )$.
\end{lemma}
\proof For $n=3$ or $n=4$ the conclusion is trivial, so assume $n
\ge 5$. Since $f(m)$ is defined only for $2 \le m \le n-1$, it is
enough to show that
\begin{equation} \label{fdiff}
f(l) - f(l-1) = f(l-1) - f(l-2) \mbox{ for all } 4 \le l \le n-1
\end{equation}
Let $i$ be the integer part of $l/2$ and $j = l - i$. Then $i \ge
2$ and either $j=i$ or $j = i + 1$, so $ j \ge 2$ too. Write
$EQ(n_1,n_2,n_3)$ for the equation (\ref{eqnf}) determined by a
particular choice of $(n_1,n_2,n_3)$. Keeping in mind that $l = i
+ j$, we have
\begin{equation} \label{EQ1} EQ( i - 1, j-1, n- l +
2): ~~~ f(l-2) + f( n - i + 1 ) + f(n - j+1 ) = C
\end{equation}
\begin{equation} \label{EQ2} EQ( i - 1, j, n- l + 1): ~~~ f(l-1) +
f( n - i + 1 ) + f(n - j ) = C \end{equation} \begin{equation}
\label{EQ3} EQ( i , j- 1, n- l + 1): ~~~ f(l-1) + f( n - i ) + f(n
- j + 1) = C
\end{equation}
\begin{equation} \label{EQ4}
EQ( i , j, n- l ): ~~~ f(l) + f( n - i  ) + f(n - j ) = C
\end{equation}
Subtract (\ref{EQ1}) from (\ref{EQ2}) to obtain
\begin{equation} f(l-1) - f(l-2)  = f( n - j + 1 ) - f( n-j )
\end{equation}
and subtract \re{EQ3} from \re{EQ4} to obtain
\begin{equation} f(l) - f(l-1)  = f( n - j + 1 ) - f( n-j )
\end{equation}
and Lemma \ref{linear f lemma} follows.
\endpf

We can now finish the proof of Theorem \ref{new}. Fix $n \ge 4$,
let $(w_j,1 \le j \le n-1)$ be a sequence of positive weights with
$w_1 = 1$.

Suppose first as in condition (i) of Theorem \ref{new} that
$(\Pi_k, 1 \le k \le n)$ is a $\Pnn$-valued fragmentation process
defined by the recursive Gibbs splitting rule derived from these
weights, with the linear selection rule, and that the distribution
of $\Pi_k$ is $p_{n,k}$ for every $k$. Then condition (ii) of
Lemma \ref{mainlem} holds for all $3 \le k \le n-1$, and in
particular for $k=3$. Lemma \ref{linear f lemma} now forces
\re{fmbc} for some $b$ and $c$, hence $w_j = w_j^{b,c}$ by Lemma
\ref{uniquew}.

Conversely, suppose that $(\Pi_k, 1 \le k \le n)$ is a
$\Pnn$-valued fragmentation process defined by the recursive Gibbs
splitting rule with the weights $w_j = w_j^{b,c}$, and the linear
selection rule. Lemma \ref{uniquew} implies that \re{fmbc} holds,
so
it is clear that condition (ii) of Lemma \ref{mainlem} holds for
$3 \le k \le n-1$ with
\begin{equation} \label{gnk}
g(n,k ) = (k-1) ( k c + n b ).
\end{equation}
Consider the inductive hypothesis that $\Pi_{k-1}$ has the
microcanonical Gibbs distribution $p_{n,k-1}^{b,c}$ with these
weights $(w_j^{b,c})$. This is true for $k=3$ by assumption.
Assuming it true for some $k$, Lemma \ref{mainlem} provides the
inductive step from $k$ to $k+1$. Thus the distribution of $\Pi_k$
is $p_{n,k}^{b,c}$ for every $2 \le k \le n-1$. Thus condition (i)
of Theorem \ref{new} is satisfied by the weights $w_j =
w_j^{b,c}$.

Condition (iii) of Theorem \ref{new}, that the reversed process is
an affine coalescent, is now read from the last sentence of Lemma
\ref{mainlem}.

\subsection{Proof of Corollary \ref{forests}}

Recall first that the distribution of an unconditioned
Galton-Watson tree, restricted to finite trees, is given by the
formula \eq \label{treeprob} \PR(T = t) = \pi(t):= \prod_{v \in
V(t)} p_{n(v,t)} \en where
\begin{itemize}
\item $t$ denotes a generic plane tree with a finite number of
nodes $\#t$; \item $V(t)$ is the set of nodes of $t$; \item
$n(v,t)$ is the number of children of node $v$ of $t$; \item $p_n$
is the probability that a node has $n$ children;
\end{itemize}
The nodes of $t$ are regarded as unlabelled. But the tree has a
root node, and the  children of each node are assigned a total
order, say from left to right. So the nodes of $t$ can be
identified or listed by some arbitrary convention, such as depth
first search, and any such convention can be used to rigorously
identify the set of nodes $V(t)$ as a subset of some ambient
countable set. See \cite{enumeration of trees,saintflour} for
background. Fix $n\ge 4$. By definition, $F_1$ is $T$ conditioned
on $\#T=n$, so
\begin{equation}
\label{tree} \PR(F_1 = t ) = \pi(t)\indic{\#t = n}/q(n)
\end{equation}
where $q(n)$ is by definition the probability that $T$ has $n$
nodes:
\begin{equation}
 \label{qn} q(n):= \PR(\#T = n) = \sum_t \pi (t)\indic{\#t = n}  = p_0^n \frac{ w_n
 }{n!}.
\end{equation}
In the last formula, read from \re{tprogbell}, the weight sequence
$(w_n)$ with $w_1 = 1$ is determined as in (\ref{g given h}) by
its exponential generating function $\sum_n w_n z^n/n!$ which is
the compositional inverse of $z p_0/g(z)$ for $g(z) :=
\sum_{n=0}^\infty p_j z^j$ the offspring generating function. The
probability \re{qn} is strictly positive for every $n \ge 1$, by
the simplifying assumption (\ref{simplify}) on  the offspring
distribution.

\newcommand{\hatF}{\widehat{F}}
Let $\hatF_2$ be the plane forest of two trees obtained by
splitting $F_1$ by deletion of a uniformly chosen random edge of
$F_1$, with subtree containing the root put to the left, and the
remaining fringe subtree put to the right. Then the distribution
of $\hatF_2$ is given by the following formula: for a generic pair
of plane trees $(t_1,t_2)$ \eq \label{f2hat} \PR( \hatF_2 =
(t_1,t_2)) = \frac{ \pi(t_1) \pi(t_2) \Sigma(t_1) } { q(n)  (n-1)
} \indic{ \#t_1 + \#t_2 = n } \en where \eq \label{sigmat}
\Sigma(t) := \sum_{v \in V(t) } r_{n(v,t)} \en with
$$
r_m:= (m+1) p_{m+1}/p_m ~~~~~~(0 \le m \le n-2)
$$
and the particular offspring distribution display in \re{off} is
characterized by the formula
\begin{equation}
\label{bcdist} r_m = b - (m-1) c ~~~~~~(0 \le m \le n-2).
\end{equation}
Formula \re{f2hat} is obtained by conditioning on which vertex $v$
of $t_1$ is the one to which $t_2$ is attached in $t$, and given
that $v$ has $m+1$ children in $t$, which  of these $m+1$ children
is the root of $t_1$. Tossing a fair coin to decide the order of
trees in $\hatF_2$ then yields $F_2$ with distribution
\begin{equation} \label{f2}
\PR( F_2 = (t_1,t_2)) = \frac{ \pi(t_1) \pi(t_2) (\Sigma(t_1) +
\Sigma(t_2) ) } {  2 q(n) (n-1) } \indic{\#t_1 + \#t_2 = n}.
\end{equation}
On the other hand, the distribution of $F_2^*$ defined by two
independent copies of $T$ conditioned to have a total of $n$ nodes
is given by
\begin{equation}
\label{f2star} \PR( F_2 ^* = (t_1,t_2)) = \frac{ \pi(t_1)
\pi(t_2)} { q_2(n)  }
 \indic{ \#t_1 + \#t_2 = n }
\end{equation} where
$$
q_2(n)= \sum_{m=1}^{n-1} q(m) q(n-m) = p_0^n \frac{ 2 B_{n,2} }{ n
! }
$$
gives the distribution  of the total progeny of the branching
process started with two individuals, as indicated in
\re{tprogbell}, with $B_{n,2} = B_{n,2}(1,w_2, \ldots, w_{n-1})$.
Condition (ii) of Corollary \ref{forests}, is the equality in
distribution \eq \label{fdist} F_2 \ed F_2^*. \en It is clear from
\re{f2} and \re{f2star} that this equality in distribution is
equivalent to the identity \eq \label{sigtree} \Sigma(t_1) +
\Sigma(t_2) = \frac{(n- 1) w_n }{B_{n,2}} \en for all pairs of
trees $(t_1,t_2)$ with $\pi(t_1)\pi(t_2) 1( \#t_1 + \#t_2 = n) >
0$, where $\Sigma(t) := \sum_{v \in V(t) } r_{n(v,t)}$ for $r_m:=
(m+1)p_{m+1}/p_m$.

\noindent (i) $\Rightarrow$ (ii). If (i) holds then $r_m = b + c -
m c$ and hence \eq \label{sigtree1} \Sigma(t_1) + \Sigma(t_2) = n
b + n c - (n-2) c = n b + 2 c \en because in every forest of two
trees with $n$ nodes the sum of the numbers of children of all
nodes is the total number of edges, which is $n-2$, and
\re{sigtree} is now read from \re{wlem}.


\noindent (ii) $\Rightarrow$ (iii) and (iv). Assuming (ii), it
follows immediately from \re{f2}, \re{f2star} and \re{fdist} that
for $k = 2$ the distribution of $\Pi_k$ generated by random
labelling of tree components of $F_k$ has the Gibbs distribution
$p_{n,k}$ with whatever weight sequence $(w_j^{b,c})$ is
associated with the distribution of the total progeny of the
branching process, and that conditionally given $\Pi_k$ the $k$
plane trees associated with these components are distributed like
independent copies of $T$ conditioned to have the sizes dictated
by the block sizes of $\Pi_k$. Suppose inductively that this is so
for some $k \ge 2$. The process of random edge deletion induces
the linear selection rule for components to split, and given that
a tree component is split, the inductive hypothesis and the
assumption for $k=2$ implies that the component is split into two
independent copies of $T$ conditioned to have the right size. The
implication (iii) $\Rightarrow$ (i) of Theorem \ref{new} now
provides the inductive step.

\noindent (iii) $\Rightarrow$ (ii) is trivial.

\noindent (ii) $\Rightarrow$ (i). This follows easily from the
identity (\ref{sigtree}) and the following Lemma:

\begin{lemma}
\label{analem} Fix $ n \ge 3$. Let $r(m)$ be a real-valued
function with domain $S= \{0,1, \ldots j \}$ for some $1 \le j \le
n-2$, such that
$$
\sum_{i=1}^ n  r(n_i) = C
$$
for some constant $C$ for each choice of $(n_1, \ldots , n_n)$
with $n_i \in S$ for all $1 \le i \le n$ and $\sum_{i = 1}^{n} n_i
= n - 2$. Then there exist real $a$ and $b$ such that $r(m) = a m
+ b$ for all $m \in S$.
\end{lemma}
\proof Consider for each $1 \le m \le j-1$ the sequence $(n_1,
\ldots , n_n)$ with the first $n-m-2$ terms equal to $1$, the next
term equal to $m$, and the last $m+1$ terms equal to $0$. This
sequence gives
$$
(n-m -2) r(1) + r(m) + (m +1) r(0) = C
$$
and the same holds for $m+1$ instead of $m$. The difference of
these two identities gives
$$
r(m+1) - r(m) = r(1) - r(0)
$$
and the conclusion follows.
\endpf

\newcommand{\Pith}{ \tilde{\Pi}_\theta}
\newcommand{\PiTh}{ \tilde{\Pi}_\Theta}
\newcommand{\Nth}{ N_\theta }
\newcommand{\tht}{ \theta/2}
\newcommand{\Fth}{ {\cal F} _ theta }
\newcommand{\Piph}{ \tilde{\Pi}_\phi}
\newcommand{\rfac}[2] { [ #1] _ {#2} }

\section{Gibbs fragmentations for random permutations in continuous time}
\label{continuous}

Given a symmetric non-negative {\em collision rate function}
$K_{i,j}$ defined for positive integers $i$ and $j$, call the
$\Pnn$-valued continuous time parameter Markovian coalescent
process $(\Pi_t, t \ge 0)$, in which each pair of clusters of
sizes $i$ and $j$ is merging at rate $K_{i,j}$, the {\em
Marcus-Lushnikov coalescent with collision kernel $K_{i,j}$}. See
\cite{aldous} for background. It is assumed throughout this
section, in keeping with the definition of a coalescent process
given in the previous section, that such a coalescent process is
started with the {\em monodisperse initial condition}. That is to
say $\Pi_0$ is the partition of $[n]$ into $n$ singletons. Both
Marcus and Lushnikov worked with the corresponding $\Pn$-valued
process rather than a $\Pnn$-valued process, but there is no
difficulty in translating results from one state-space to the
other, by application of the standard criterion for a function of
a Markov process to be Markov. Lushnikov \cite{lushnikov} found
the remarkable result that for a collision kernel of the form
$K_{i,j} = i f(j) + j f(i)$ for each $t > 0$ the distribution of
$\Pi_t$ is of the form
\begin{equation}
\PR( \Pi_t = \pi ) = \sum_{k=1}^n q_{n,k}(t) p_{n,k}(\pi ; w_j(t),
j = 1, 2, \ldots )
\end{equation} where $p _{n,k}(\pi ; w_j, j = 1, 2, \ldots )$
denotes the microcanonical Gibbs distribution on $\Pnk$ with
weights $w_j$, and the functions $q_{n,k}(t) = \PR(\#\Pi_t = k )$
and the weights $w_j(t)$ are determined by a system of
differential equations. As mentioned earlier, Hendriks et al.
\cite{hses85} showed that for $K_{i,j} = a + b (i + j )$ for
constants $a$ and $b$ the $w_j(t)$ can be chosen independently of
$t$ as $w_j(t) = w_j$ where $w_j$ is determined by $a$ and $b$ via
formula (\ref{wjbc}) for $c = a/2$.

Kingman \cite{ki82co} studied the particular case of the
Marcus-Lushnikov coalescent with $a = 1$ and $b = 0$. In this
process, at any given time $t$, given that $\#\Pi_t = k$, each of
the $k(k-1)/2$ cluster pairs in existence at time $t$ is merging
at rate $1$. Call this process with state space $\Pnn$ {\em
Kingman's $n$-coalescent}, Motivated by applications to genetics,
Kingman \cite{ki82co} proposed the following construction. Given a
coalescent process $(\Pi_t, t \ge 0 )$, suppose that each cluster
of $\Pi_t$ is subject to mutation at rate $\theta/2$ for some
$\theta >0$. Now define a random partition $\Pith$ of $[n]$ by
declaring that $i$ and $j$ are in the same block of $\Pith$ if and
only if no mutation affects the clusters containing $i$ and $j$ in
the interval $(0, \tau_{ij})$ where $\tau_{ij}$ is the {\em
collision time} of $i$ and $j$ in the coalescent process $(\Pi_t,
t \ge 0 )$, that is the first time $t$ that $i$ and $j$ are in the
same cluster of $\Pi_t$. Kingman obtained the following result:

\begin{prp}{king}
{\em (Kingman \cite{kingman82esf})} Suppose that $(\Pi_t, t \ge 0
)$ is Kingman's $n$-coalescent. Then \begin{equation} \label{pith}
\PR( \Pith = \pi ) = { \theta ^{k-1} \over \rfac{\theta + 1 }{n-1}
} \prod_{i=1}^k (n_i - 1)! \end{equation} for each partition $\pi$
of $[n]$ into $k$ components of sizes $n_1, \ldots, n_k$
\end{prp}

The distribution of $\Pith$ defined by (\ref{pith}) first appears
in \cite{ewens} and is known as Ewens' Sampling Formula with
parameter $\theta$. This distribution has long been recognized as
an essential tool in population genetics (see, e.g.
\cite{durrett-schw} recently), and has been applied in a wide
variety of contexts in probability. Note that this distribution is
a particular mixture over $k$, with mixing coefficients depending
on $\theta$, of the microcanonical Gibbs distributions on $\Pnk$
with weights $(j-1)!$, as interpreted in Example \ref{exm2}.
Recall that this distribution has been constructed starting from
$(\Pi_t,t\ge 0)$, where due to Example \ref{exm3}, for all $t$,
$\Pi_t$ is a mixture over $k$, with mixing coefficients depending
on $t$, of the microcanonical Gibbs distributions on $\Pnk, 1 \le
k \le n$, with the different weight sequence $(j!, j \ge 1)$. It
does not seem obvious from a combinatorial perspective why there
should be such a connection between the Gibbs models with these
two weight sequences.


The random partition $\Pith$ (which is sometimes referred to as
the random allelic partition), and more generally the
fragmentation process $(\Pith, \theta \ge 0)$ discussed below, can
be defined starting from any coalescent $(\Pi_t)$, but there seems
to be a manageable formula for the distribution of $\Pith$ only
for Kingman's coalescent. See however the recent work of M\"ohle
\cite{mohle} where an explicit recursion is given for the random
allelic partition obtained from a $\Lambda$-coalescent (i.e.,
coalescent with multiple collisions). See also the related work of
\cite{dong} as well as \cite{bbs1, bbs2} which has some explicit
asymptotic formulae in the particular case of a beta-coalescent.

\medskip As a development of Kingman's result, there is the
following proposition. See also \cite{gnedin-pitman-esf} (or
\cite[Exercise 5.2.1]{saintflour}) for an alternative
construction.

\begin{proposition}\label{coalth} There exists a Gibbs
fragmentation process $(\Pith, \theta \ge 0)$ with weight sequence
$((j-1)!, j \ge 1)$ such that for each $\theta > 0$ the
distribution of $\Pith$ is the Gibbs distribution on $\Pnn$ with
these weights as displayed in \emph{(\ref{pith})}.
\end{proposition} \proof Given the path of a Kingman coalescent
process $(\Pi_t, t \ge 0)$, construct a random tree $\TT$ as
follows. Let the vertices of the tree $\TT$ be labelled by the
random collection $\VV$ of all subsets of $[n]$ which appear as
clusters in the coalescent at some time in its evolution. Because
the coalescent develops via binary mergers, starting with $n$
singletons and terminating $\Pi_t = [n]$ for all sufficiently
large $t$, the set $\VV$ comprises the collection of all $n$
singleton subsets of $[n]$, which are the {\em leaves} of the
tree, the whole set $[n]$ which is the {\em root} of the tree, and
$n-2$ further subsets of $[n]$, whose identities depend on how the
coalescent evolves, which are the {\em internal vertices} of the
tree. The tree $\TT$ has $ n + 1 + (n-2) = 2 n - 1$ vertices all
together. Associate with each subset $v$ of $[n]$ that is a vertex
of the tree the time $t(v)$ at which the coalescent forms the
cluster $v$. Thus $t(v) = 0$ if and only if $v$ is one of the $n$
singleton leaf vertices, $t([n]) = \inf \{t : \# \Pi_t = 1 \}$,
and the collection of times $t(v)$ as $v$ ranges over the $n-1$
non-leaf vertices of the tree is the set of times $t$ at which the
process $( \# \Pi_t , t \ge 0 )$ experiences a downward jump. For
each non-leaf vertex $v$ in $\TT$, let there be exactly two edges
of $\TT$ directed from $v$ to $v_1$ and $v_2$, where $v_1$ and
$v_2$ are the two clusters which merged to form $v$. Let each
vertex $v$ of $\TT$ be placed at height $t(v)$ equal to the time
of its formation, and for $i = 1,2$ regard the directed edge from
$v$ to $v_i$ as a segment of length $t(v) - t(v_i)$. Now Kingman's
construction of $\Pith$ amounts to supposing that there is a
Poisson process of cut points on the edges of this tree, with rate
$\theta/2$ per unit length, and identifying the blocks of $\Pith$
with the restrictions to the set of $n$ leaves of $\TT$
(identified with $[n]$) of the components of the random forest
obtained by cutting segments of $\TT$ at the Poisson cut points.
Now conditionally given the tree $\TT$, construct the Poisson cut
points simultaneously for all $\theta > 0$ so that for each edge
of the tree of of length $\ell$ the moments of cuts of that edge
form a homogeneous Poisson process of rate $\ell \theta/2$, and
these processes are independent for different edges. Then $\Pith$
has been constructed simultaneously for each $\theta >0$ in such a
way that $\Pith$ is obviously a refinement of $\Piph$ for $\theta
> \phi$. Since with probability one there are no ties between the
times of cuts on different segments, it is clear that the process
$(\Pith, \theta \ge 0)$ develops by binary splits. Thus $(\Pith,
\theta \ge 0)$ is a Gibbs fragmentation process. \endpf

While the one-dimensional distributions of this process $(\Pith,
\theta \ge 0)$ are given by Ewens' sampling formula (\ref{pith}),
the two and higher dimensional distributions seem difficult to
describe explicitly. In particular, a calculation of the simplest
transition rate associated with the process $(\Pith, \theta \ge
0)$, provided below, shows that this rate depends on $\theta$. It
seems quite difficult to give a full account of all transition
rates of $(\Pith, \theta \ge 0)$, though their general form can be
described and a method for their computation for small $n$ will be
indicated. For $n \ge 2$ the process $(\Pith, \theta \ge 0)$ turns
out to be non-Markovian, so its distribution is not determined by
its transition rates.

In connection with Proposition \ref{coalth} and such calculations,
the following problem arises:

\begin{problem}\label{pbmfacwts}
Does there exist for each $n$ a $\Pnn$-valued Gibbs fragmentation
process $( \Pi_k , 1 \le k \le n )$ with weight sequence $((j-1)!,
j \ge 1 )$?
\end{problem}


\subsection{Calculations with the tree derived from Kingman's coalescent.}

The following calculations (Proposition \ref{calc}) show that for
$n\ge 4$, the discrete-time chain embedded in $(\Pith, \theta \ge
0)$ (that is, the sequence of successive states of
$(\Pith,\theta>0)$, or its discrete skeleton) does not provide a
solution to Problem \ref{pbmfacwts}.

Let $\TTn$ denote the random tree derived as in the proof of
Proposition \ref{coalth} from Kingman's $n$-coalescent $(\Pi_t, t
\ge 0 )$, and recall the definition of the fragmentation process
$(\Pith,\theta>0)$. Let $\Theta$ be the time of the first cut in
this process, and let $\PiTh$ be the state of the fragmentation at
this random time. Thus almost surely $\PiTh$ is a partition with
two blocks.

\begin{proposition} \label{calc}
The law of $\Theta$ is determined by
\begin{equation} \label{density Theta}
\PR( \Theta \in  d \theta ) /d \theta = {(n-1) ! \over
\rfac{\theta + 1 }{n-1} }\, \sum_{i = 1}^{n-1} { 1 \over i +
\theta}
\end{equation}
For $\pi$ with two components of sizes $n_1$ and $n_2$
\begin{equation} \label{PiTh giv Theta}
\PR( \PiTh = \pi \giv \Theta = \theta) ={\sum_{i=1}^{n-1} (i +
\theta)^{-1} { n - 1 \choose i }^{-1} \left[  { n_2 - 1\choose i -
1} + { n_1 - 1 \choose i - 1} \right] \over { n \choose n_1 }
\sum_{j = 1}^{n-1} (j + \theta)^{-1} }
\end{equation}
and
\begin{equation} \label{PiTh distr} \PR( \PiTh = \pi ) = { (n-1)!  \over { n \choose
n_1 } } \, \int_0^\infty { d \theta \over \rfac{ \theta +1 }{n-1}
} \, \sum_{i=1}^{n-1} { \left[  { n_1-1 \choose i - 1} + { n_2-1
\choose i - 1} \right] \over { n - 1  \choose i } (i + \theta) }
\end{equation}
\end{proposition}

\begin{proof} We only provide a sketch of the calculations leading to this result as they are somewhat tedious.
For $1 \le k \le n$ let $T_k = \inf \{t : \# \Pi_t = k \}$. Let
$S_i=(i+1)(T_i-T_{i+1})$, which we call the $i\th$ stratum of the
tree. It is clear that the total length of all segments in the
tree $\TT$ is
\begin{equation} \label{L_n}
L_n := \sum_{i=1}^{n-1} (i+1) (T_{i} - T_{i+1})
\end{equation}
From the definition of the underlying coalescent process $(\Pi_t,
t \ge 0 )$, the random variable $T_{i+1} - T_i$ has exponential
distribution with rate $i(i+1)/2$, and these random variables are
independent for $1 \le i \le n-1$. It follows that
\begin{equation} E \exp ( - { \theta
\over 2 } L_n) = \prod_{i=1}^{n-1} { i(i + 1) / 2 \over
\theta(i+1)/2 + i (i+1 )/2 } = {(n-1) ! \over \rfac{\theta + 1
}{n-1} } \label{lt}
\end{equation}
where $\rfac{\theta + 1 }{n-1} = \prod_{i=1}^{n-1} (\theta + i )
$. On the other hand, given $L_n$, the Poisson process with rate
$\theta/2$ per unit segment length in the tree has no points with
probability $\exp ( - (\theta /2) L_n)$. So the expectation
calculated in (\ref{lt}) is just the probability that $\Pith$ is
the partition of $[n]$ with one component, or in other words that
$\Theta >\theta$. Thus (\ref{density Theta}) follows by
differentiation.

Now, let $I$ denote the index of the stratum in which the first
cut point falls at time $\Theta$. Then it follows from the
representation  of $L_n$ as the sum of independent exponential
variables $L_n = \sum_{i=1}^{n-1}S_i$ that the sum over $i$ in
(\ref{density Theta}) corresponds to summing over the possible
values $i$ of $I$. That is, for $1 \le i \le n-1$,
\begin{equation} \PR( \Theta \in  \theta , I = i ) /d \theta =
{(n-1) ! \over \rfac{\theta + 1 }{n-1} }\, { 1 \over i + \theta}
\end{equation}
and hence
\begin{equation} \label{igivth} \PR( I =
i \giv \Theta = \theta) = { (i + \theta)^{-1 } \over \sum_{j =
1}^{n-1} { (j + \theta)^{-1}} }.
\end{equation}
Observe now that given $\Theta = \theta$ and $I = i$, the
partition $\PiTh$ consists of two components, obtained as the
restriction to $[n]$, identified as the set of leaves of the tree
$\TT$, of the two components of $\TT$ separated by the cut at time
$\Theta$ in stratum $i$ of $\TT$. To be precise, $\PiTh = \{C, [n]
- C \}$ where $C$ is the cluster of $\Pi_t$ in existence during
the time interval $(T_{i+1}, T_i)$ corresponding to the segment of
$\TT$ which is cut at time $\Theta$. This $C$ is one of the
clusters of $\Pisk{i+1}$, where $(\Pisk{k}, 1 \le k \le n)$  is
the discrete skeleton of $(\Pi_t, t \ge 0)$. Now by construction
of the Poisson cutting process, and the fact that the discrete
skeleton $(\Pisk{k}, 1 \le k \le n)$ of $(\Pi_t, t \ge 0)$ is
independent of $(\# \Pi_t, t \ge 0)$, it is clear that the
conditional distribution of $\Pisk{i+1}$ given $\Theta = \theta$
and $I = i$ is identical to its unconditional distribution, that
is the Gibbs distribution on ${\cal P}_{[n,i+1]}$ with weights
$(j!, j \ge 1)$, and that $C$ is one of the $i+1$ components of
$\Pisk{i+1}$ picked by a mechanism independent of the sizes of
these components. Therefore,
\begin{equation}
\PR( \# C = n_1 \giv
\Theta = \theta, I = i ) = \PR( \# C_{i+1} = n_1 )
\end{equation}
where $C_{i+1}$ is a random component of $\Pisk{i+1}$. After some
combinatorics, it follows easily that the conditional distribution
of $\PiTh$ given $\Theta = \theta$ and $I = i$ is given by
\begin{equation}
\PR( \PiTh = \pi \giv \Theta = \theta, I = i ) = { n  \choose n_1
} ^{-1} { n - 1 \choose i }^{-1} \left[ { n_2 - 1 \choose i - 1} +
{ n_1 - 1 \choose i - 1} \right].
\end{equation}
Combining this expression with (\ref{igivth}) shows that the
conditional distribution of $\PiTh$ given $\Theta = \theta$ is
given by (\ref{PiTh giv Theta}). We now easily obtain (\ref{PiTh
distr}) from (\ref{PiTh giv Theta}) and (\ref{density Theta}) by
integration. \end{proof}

We now briefly explain why it can be deduced from the explicit
formula (\ref{PiTh distr}) that the discrete chain embedded in
$(\Pith,\theta>0)$ is not a Gibbs fragmentation. We must simply
show that (\ref{PiTh distr}) does not coincide with the Gibbs
microcanonical distribution $p_{n,2}$ associated with the weight
sequence $w_j=(j-1)!$. Let $J_\Theta$ denote the size of a
component of $\PiTh$ picked by the toss of a fair coin independent
of $\PiTh$. Then using the above, if $n_1=1$ and $n_2=n-1$,
$$
\PR(J_\Theta =1)= \frac12 (n-2)!
\int_0^{\infty}\frac{d\theta}{[\theta+1]_{n-1}}\left(\sum_{i=1}^{n-1}
\frac{i}{i+\theta} + \frac1{1+\theta}\right)
$$
After a few lines of algebra, using some integration by parts and
some partial fractions, we can conclude that
\begin{equation}\label{JTh}
\PR(J_\Theta=1)=\frac12(n-1)!\sum_{i=1}^{n-1}a_{n,i} \log i
+\frac12 .
\end{equation}
where $a_i=(-1)^{i-1} {n\choose i-1}/ n!$.

On the other hand, let $\Pi$ have the Gibbs $(n,2,w)$ distribution
with $w_j=(j-1)!$, and $J$ is the size of a randomly picked
component. Remark that by decomposing on the size of the cycle
containing 1, $$ B_{n,2} =\sum_{j=1}^{n-1} {n-1 \choose
j-1}(j-1)!(n-1-j)! =(n-1)!H_{n-1}$$ where $H_{n-1} := \sum_{j =
1}^{n-1}1/j$. Since the number of permutations with exactly two
cycles one of which has size 1 is $n(n-2)!$, we conclude that $
\PR(J=1)= \frac12\frac{n}{(n-1)H_{n-1}} $. This is incompatible
with (\ref{JTh}). Indeed if this was to be equal to right-hand
side in (\ref{JTh}) one would get $\log (\prod_{i=1}^{n-1}
i^{r_i}) =q$ for some rational number $q$ and $r_i=(-1)^i
{n\choose i-1} \in \mathbf{Z}$, and thus $e^q=q'$ for some (other)
rational numbers $q$ and $q'$. This contradicts the transcendence
of $e$.

Thus the distribution of the partition of $[n]$ into two parts
obtained at the time of the first split is not the common
distribution of $\Pith$ given $\# \Pith = 2$ for all $\theta >0$.
In particular, for $n=4$ the formulas above give: $$ \PR( J_\Theta
= 1 ) = \PR( J_\Theta = 3 ) = 3(- \log 2 + \frac12 \log 3)+\frac12
$$ and $$
 \PR( J_\Theta = 2 ) = 6\log 2 -3 \log 3
$$
whereas
$$
\PR(J=1) = \PR(J=3) = 4/11
$$
$$
\PR(J=2) =3/11 
$$ We conclude that the discrete skeleton of $(\Pith)_{\theta \ge
0}$, i.e., the discrete fragmentation chain embedded in it, does
not give a discrete Gibbs fragmentation associated with
$w_j=(j-1)!$.

\subsection{A reformulation with walks on the symmetric group}

In view of the combinatorial interpretation of Example \ref{exm2},
Problem \ref{pbmfacwts} can be restated as:

\begin{problem}\label{permut pbm}
Does there exist for each $n$ a sequence of random permutations $(
\sigma_k , 1 \le k \le n )$ such that $\sigma_k$ has uniform
distribution on the set of permutations of $[n]$ with $k$ cycles,
and for $k \le \ell$ the partition generated by the cycles of
$\sigma_\ell$ is a refinement of $\sigma_k$?
\end{problem}

This problem may be partially reformulated in terms of random
walks on the symmetric group. Suppose we consider the Cayley graph
$G_n$ of the symmetric group induced by the set of generators
$S=\{\text{all transpositions}\}$, that is, we put an edge between
two permutations $\sigma$ and $\pi$ if and only if $\sigma$ may be
written as $\sigma=\tau\cdot \pi$ for some transposition $\tau$.
It is well-known that multiplying a permutation by a transposition
can only result in a coagulation or a fragmentation in the cycle
structure. More precisely, suppose $C=(x_1,\ldots,x_k)$ is a cycle
of the permutation $\pi$. If we multiply by the transposition
$\tau=(x_i,x_j)$ then the resulting cycle structure in $\sigma$ is
the same as that of $\pi$ except that $C$ breaks into
$(x_1,\ldots,x_{i-1},x_j,x_{j+1},\ldots,x_k)$ on the one hand and
$(x_i,\ldots,x_{j-1})$ on the other hand. Conversely, suppose
$C=(x_1,\ldots,x_k)$ and $C'=(y_1,\ldots,y_l)$ are two cycles of
$\pi$ and we multiply $\pi$ by the transposition $(x_i,y_j)$. In
the resulting permutation, $C$ and $C'$ will be replaced by a
unique cycle $C''=(x_1,\ldots,
x_{i-1},y_j,y_{j+1},\ldots,y_{j-1},x_i,\ldots,x_k)$. In
particular, any (random) walk on $G_n$ may be viewed as a
coagulation and fragmentation process on the cycle structure of
the permutation. Moreover, a well-known result due to Cayley
states that if $\sigma$ is a permutation then the \emph{graph
distance} between $\sigma$ and the identity permutation $I$ (i.e.,
the minimum number of edges one must cross to go from $I$ to
$\sigma$ on $G_n$) is simply $n-\#\text{cycles}$ of $\sigma$. As a
consequence, by considering a time-reversal of the process, to
solve Problem \ref{permut pbm}, it is enough to construct a random
process $(\sigma_k)_{0\le k \le n-1}$ on $G_n$ which has the
following two properties:
\begin{enumerate} \item The sequence $(\sigma_k,0\le k\le n-1)$ is a random walk on $G_n$, in
the sense that if $\sigma_k=\sigma$ at the next stage $\sigma$ can
only jump to one the neighbors of $\sigma$.

\item The permutation $\sigma_k$ has the following marginal
distribution: at stage $k$ the distribution of $\sigma_k$ is
uniform on the sphere of radius $k$ about the identity, that is
the set of all permutations whose distance to the identity is $k$.
\end{enumerate}

\noindent Property 1 ensures that the cycles of $\sigma_k$ perform
a coagulation-fragmentation process. In conjunction with property
2, since $\sigma_k$ must be at distance $k$, it must be the case
that all jumps of $\sigma$ are produced by some fragmentation.
Moreover, by Cayley's result for the distance of a permutation to
the identity, if $\sigma_k$ is uniform on the sphere of radius $k$
then its cycle structure is a realization of the Gibbs
distribution (\ref{gibbsnk}) with weight sequence $w_j=(j-1)!$
(Note however this is not strictly equivalent to Problem
\ref{permut pbm} since not all fragmentations at the level of
partitions can be represented by moving along some edge of $G_n$.
For instance, it is impossible to get from the permutation $(1\ 2\
3\ 4)$ to the permutation $(1\ 3)(2\ 4)$ in one step).

\mn In this context, a very natural process to consider is the
simple random walk on $G_n$, conditioned to never backtrack. In
other words, starting from the identity, at each step choose
uniformly among all edges that lead from distance $k$ to distance
$k+1$. Although this may seem a very natural candidate for
properties 1 and 2, this is far from being the case. Much is known
about this process, and in particular it has been shown in
\cite{nb} that the distribution of this process at time $k=\lfloor
an\rfloor$ for any $0<a<1$ is asymptotically singular with respect
to the uniform distribution on the sphere of radius $k$.

\mn \textbf{Acknowledgements}

\mn We thank Jomy Alappattu for carefully reading a draft of this
paper and for pointing out some mistakes, as well as some helpful
discussions. We thank the referees for some useful suggestions.


\end{document}